\documentclass[11pt,a4paper,reqno]{amsart}
\usepackage[T1]{fontenc}
\usepackage{lmodern}
\usepackage{amsmath,amssymb,amsthm,mathtools,mathrsfs}
\usepackage{needspace}
\usepackage[a4paper,left=27mm,right=27mm,top=27mm,bottom=28mm,headheight=14pt,headsep=8mm,footskip=12mm]{geometry}
\usepackage{microtype}
\usepackage{xcolor}
\usepackage[colorlinks,linkcolor=blue, citecolor=red,anchorcolor=blue, pagebackref=true ]{hyperref}
\usepackage{bookmark}

\newcommand{\Fp}{\mathbb F_p}
\newcommand{\PG}{\textup{PG}}
\newcommand{\AG}{\textup{AG}}
\theoremstyle{plain}
\newtheorem{theorem}{Theorem}[section]
\newtheorem{lemma}[theorem]{Lemma}
\newtheorem{proposition}[theorem]{Proposition}
\newtheorem{corollary}[theorem]{Corollary}
\theoremstyle{definition}
\newtheorem{definition}[theorem]{Definition}
\theoremstyle{remark}
\newtheorem{remark}[theorem]{Remark}
\newtheorem*{projectiveform}{Theorem \ref{thm:main}$^\ast$}
\numberwithin{equation}{section}

\allowdisplaybreaks[1]
\title[The cylinder conjecture and R\'edei's conjecture]{The strong cylinder conjecture and R\'edei's factorization conjecture in three dimensions}
\author{Haihua Deng}
\author{Weicong Li}
\author{Weiye Zhang$^\ast$}
\thanks{$^\ast$Corresponding author}
\address[Haihua Deng and Weiye Zhang]{School of Mathematical Sciences, Zhejiang University, Hangzhou 310058, China}
\email{haihua.deng@zju.edu.cn, zhangweiye@zju.edu.cn}
\address[Weicong Li]{Department of Mathematics, School of Sciences, Great Bay University, Dongguan, China}
\email{liweicong@gbu.edu.cn}
\date{}

\makeatletter
\def\ps@plain{%
  \let\@mkboth\@gobbletwo
  \let\@oddhead\@empty
  \let\@evenhead\@empty
  \def\@oddfoot{\normalfont\small\hfil\thepage\hfil}%
  \let\@evenfoot\@oddfoot
}
\makeatother
\makeatletter
\renewcommand{\subsection}{%
  \@startsection{subsection}{2}%
    \z@{.5\linespacing\@plus.7\linespacing}{.25\linespacing}%
    {\normalfont\bfseries}}
\renewcommand{\subsubsection}{%
  \@startsection{subsubsection}{3}%
    \z@{.5\linespacing\@plus.7\linespacing}{.25\linespacing}%
    {\normalfont\itshape}}
\makeatother

\begin{document}
\setcounter{tocdepth}{1}

\begin{abstract}
In this paper, we prove the strong cylinder conjecture, that is, for every prime $p$, a set of $p^2$ points in $\Fp^3$ that meets every affine plane in a multiple of $p$ points is a union of $p$ distinct parallel affine lines. We use this result to resolve R\'edei's 1970 factorization conjecture in three dimensions: if $A,B\subseteq\Fp^3$ contain zero and every vector in $\Fp^3$ has a unique expression $a+b$ with $a\in A$ and $b\in B$, then one of the two sets is contained in a plane through zero.
\end{abstract}
\maketitle
\tableofcontents
\thispagestyle{plain}
\pagestyle{plain}

\section{Introduction}\label{sec:intro}

A classic question in finite geometry is how many directions a set of affine points determines. Let $p$ be a prime. Two distinct points $x,y\in\Fp^m$ determine the direction $\Fp(y-x)$. A direction is called \emph{undetermined} by a set if no two of its points determine it. In an affine plane, $p$ collinear points determine just one direction. For odd $p$,  R\'edei and Megyesi showed that $p$ noncollinear points determine at least $(p+3)/2$ of the $p+1$ directions; see \cite[Theorem~24']{RedeiBookEnglish}. Lov\'asz and Schrijver \cite{LS} gave a simpler proof, and Somlai \cite{SomlaiDirections} later gave another proof. For more details, see the survey \cite{DBsurvey}.

R\'edei studied directions in connection with factorizations of finite abelian groups. We write $A\oplus B=G$ if every element of $G$ has a unique expression $a+b$ with $a\in A$ and $b\in B$. The translates $A+b$, for $b\in B$, then partition $G$. Such a partition is called a \emph{tiling}, and the factorization is called \emph{normalized} if $0\in A\cap B$. In a factorization $A\oplus B=\Fp^2$ with $|A|=|B|=p$, one factor is an affine line, through zero if the factorization is normalized; see \cite[Satz~4]{RedeiGaps1947} and \cite[Theorem~23]{RedeiBookEnglish}. This is the planar case of R\'edei's theorem on factors of prime size \cite{RedeiFactors}.

In $\Fp^3$, two factors of size greater than one must have sizes $p$ and $p^2$. R\'edei asked whether, in a normalized factorization of this group, one factor must lie in a plane through zero \cite{RedeiBook}. This is Problem~5 in the English edition \cite[Chapter ``Some unsolved problems'']{RedeiBookEnglish}. Szab\'o and Ward \cite{SzaboWard} verified the conjecture for $p\le11$, and Szab\'o \cite[Theorem~1]{Szabo13} extended the verification to $p=13$. Ball \cite{BallFew} suggested that direction problems in three dimensions could help resolve this factorization problem. R\'edei's result on replacing factors \cite{RedeiFactors} implies that every direction determined by the smaller factor is undetermined by the larger one; see \cite[Section~2]{SzaboWard} and \cite[Proposition~6.2.5]{MattheusThesis}.

If a set of $p^2$ points in $\Fp^3$ does not determine a direction, each of the $p^2$ affine lines in that direction contains exactly one point of the set. Taking these lines to be vertical gives a graph $z=f(x,y)$. Ball and Lavrauw \cite[Theorem~2.2]{BL} proved that if at least $p$ directions are undetermined, every affine plane meets the set in a multiple of $p$ points. A set with this plane intersection property is called \emph{$p$-divisible}. The union of $p$ distinct parallel affine lines is called a \emph{cylinder}. The \emph{strong cylinder conjecture} asks whether every $p$-divisible set of $p^2$ points in $\Fp^3$ is a cylinder \cite{BBGSW,BallFew}. The \emph{weak cylinder conjecture} asks the same for a set of $p^2$ points with at least $p$ undetermined directions \cite{Ball}.

Mattheus \cite[Conjecture~6.2.7 and Theorem~6.2.11]{MattheusThesis} gave a conditional reduction of R\'edei's problem. If every nonplanar set of $p$ points in $\Fp^3$ determines at least $p+1$ directions, then the larger factor has at least $p+1$ undetermined directions whenever the smaller factor spans $\Fp^3$. A cylinder theorem under this hypothesis would suffice, since projection along the parallel lines gives a factorization of the affine plane. See also \cite[Section~7]{KMNS}.

Several cases of the cylinder conjectures were known. For $p>3$, Sziklai \cite[Proposition~11]{SziklaiDirections} proved that a set of $p^2$ points determining fewer than $p(p+5)/2$ directions is an affine plane or a cylinder. De Beule, Demeyer, Mattheus, and Sziklai \cite{BDMS} proved the cylinder conclusion for $p\le13$ under the assumption that at least $p+1$ directions are undetermined. They also proved the strong conjecture for $p=2,3$. Kurz and Mattheus \cite[Propositions~41 and~52]{KM} corrected the proof for $p=5$ and proved the case $p=7$. Kiss, Mark\'o, Nagy, and Somlai \cite[Theorem~1.12]{KMNS} showed that the indicator function of every $p$-divisible set of $p^2$ points is an integer linear combination of indicator functions of cylinders. Bartoli and Giannoni \cite[Theorem~1.1]{BG} proved the weak conjecture for polynomial graphs of total degree at most $d$, where $d\ge3$ and $p>d^{9/2}$.

We prove the following theorem, confirming the strong cylinder conjecture \cite{BBGSW} \cite[Conjecture~5.2]{Ball}.

\begin{theorem}\label{thm:main}
Let $p$ be a prime, and let $S\subseteq\Fp^3$ satisfy
\begin{equation}\label{eq:hyp_af}
 |S|=p^2,\qquad |S\cap\Gamma|\equiv0\pmod p
\end{equation}
for every affine plane $\Gamma$ of $\Fp^3$. Then $S$ is a cylinder.
\end{theorem}

Using Theorem~\ref{thm:main}, we also prove R\'edei's factorization conjecture \cite{RedeiBook} \cite[Problem 5]{RedeiBookEnglish}.

\begin{theorem}\label{thm:redei}
Let $p$ be a prime. If $A\oplus B=\Fp^3$ is normalized, then $A$ or $B$ is contained in a plane through the zero vector.
\end{theorem}

Every elementary abelian group of order $p^3$ is isomorphic to the additive group of $\Fp^3$, in which the subgroup generated by a set is its $\Fp$-linear span. Thus, Theorem \ref{thm:redei} is equivalent to the group-theoretic formulation in \cite{SzaboWard}, i.e., in every normalized factorization, one factor generates a proper subgroup.

This paper is organized as follows. Section~\ref{sec:prelim} introduces the notation and preliminary results. In Section~\ref{sec:projection}, we study how $S$ meets lines and planes and obtain polynomial identities from these intersection counts. Section~\ref{sec:exact} uses these identities to determine the possible sizes of plane intersections and to establish the degree bounds needed later. In Section~\ref{sec:ci}, we establish polynomial identities that will be used to construct a quadratic equation in Section~\ref{sec:quadric}. Section~\ref{sec:quadric} combines these results with geometric arguments to prove Theorem~\ref{thm:main}. Finally, Section~\ref{sec:redei} applies Theorem~\ref{thm:main} to prove R\'edei's factorization conjecture.


\section{Preliminaries}\label{sec:prelim}

Throughout this paper, let $\Fp$ be the finite field with $p$ elements where $p$ is a prime. Denote by $K=\overline{\Fp}$ the algebraic closure of $\Fp$. For an extension field $\mathbb E$ of $\Fp$, set $\mathbb E^\times=\mathbb E\setminus\{0\}$. We fix the notation used below and explain the correspondence between affine and projective space.

\subsection{Affine spaces, projective spaces, and point-plane duality}\label{sec:affine-projective}

We write $\PG(m,K)$ for the projective $m$-space over $K$, whose points are the one-dimensional subspaces of $K^{m+1}$. For a finite field $\mathbb F_q$, we write $\PG(m,q)=\PG(m,\mathbb F_q)$. For a nonzero vector $v=(v_0,\ldots,v_m)\in K^{m+1}$, we write $[v]=(v_0:\cdots:v_m)$ for the corresponding projective point. For a set $X$ of projective points, $\langle X\rangle$ denotes its projective span; for a set of vectors, the same notation denotes its linear span over the ground field.

A point, line, or plane in $\PG(m,K)$ is called \emph{$\Fp$-rational} if it is defined over $\Fp$. We identify $\PG(m,p)$ with the set of $\Fp$-rational points of $\PG(m,K)$.

We next introduce the affine notions needed in the sequel, following \cite[Section~2]{KM}.

\begin{definition}
Let $V=\Fp^3$. A \emph{direction} in $V$ is a one-dimensional subspace $L\le V$. An \emph{affine line} with direction $L$ is a coset $a+L$ for some $a\in V$.


For each direction $L$, consider the natural quotient map
\[
\pi_L:V\longrightarrow V/L,\qquad
a\longmapsto a+L.
\]
The \emph{preimage} of $a+L\in V/L$ under $\pi_L$ is
\[
\pi_L^{-1}(a+L)
=\{a+\ell:\ell\in L\}=a+L.
\]
Thus, the preimages of points under $\pi_L$ are precisely the affine lines with direction $L$.
\end{definition}

We now describe the quotient map in coordinates. Fix a direction $L=\Fp v$. Since both \eqref{eq:hyp_af} and the cylinder property are invariant under affine changes of coordinates, we may choose coordinates $(x,y,z)$ such that $v=(0,0,1)$. Identifying $V/L$ with $\Fp^2$, the quotient map then takes the form
\[
\pi_L:V\longrightarrow\Fp^2,\qquad
(x,y,z)\longmapsto(x,y).
\]


Let $\widehat V=\Fp^4$, with dual space $\widehat V^*=\operatorname{Hom}_{\Fp}(\widehat V,\Fp)$. The dual projective space $\PG^*(3,p)$ consists of the one-dimensional subspaces of $\widehat V^*$. Choose homogeneous coordinates $(X:Y:Z:T)$ on $\PG(3,p)$ and let $\Pi_\infty$ be the plane $T=0$. In these coordinates, we identify $V$ with $\PG(3,p)\setminus \Pi_\infty$ via
\[
(x,y,z)\longmapsto(x:y:z:1),
\]
see \cite[Section~1]{KMNS}.

Under this identification, we regard the set $S$ in \eqref{eq:hyp_af} as a subset of $\PG(3,p)\setminus\Pi_\infty$. For each projective plane $\Pi\ne\Pi_\infty$, the corresponding affine plane is $\Pi\setminus\Pi_\infty$. Since $S\cap\Pi_\infty=\varnothing$, the plane-intersection condition in \eqref{eq:hyp_af} is equivalent to $|S\cap\Pi|\equiv0\pmod p$ for every projective plane $\Pi$.

The $p$ parallel affine lines of a cylinder extend to $p$ projective lines through a common point $v\in\Pi_\infty$; we also call their union without the common point $v$ a \emph{cylinder}. Thus, Theorem~\ref{thm:main} is equivalent to the following projective formulation.

\begin{projectiveform}
Let $p$ be a prime, and let $\Pi_\infty$ be a plane of $\PG(3,p)$. Suppose that $S\subseteq\PG(3,p) \setminus \Pi_\infty$ satisfies
\begin{equation}\label{eq:hyp}
|S|=p^2,\qquad
|S\cap\Pi|\equiv0\pmod p
\end{equation}
for every projective plane $\Pi$. Then $S$ is a cylinder.
\end{projectiveform}

Lemma~\ref{lem:infinity} also shows that the cylinder property is unchanged when $\Pi_\infty$ is replaced by any projective plane disjoint from $S$.

We use the standard point-plane duality between $\PG(3,p)$ and $\PG^*(3,p)$; see \cite[Problem~4.18 (a), (b)]{Fulton}. Each nonzero linear form $\lambda\in\widehat V^*$ determines a projective plane $\Pi=\{[w]\in\PG(3,p):\lambda(w)=0\}$, whose dual point is $\Pi^*=[\lambda]$. In homogeneous coordinates, this correspondence is given by
\[
\Pi:\ aX+bY+cZ+dT=0
\quad\longleftrightarrow\quad
\Pi^*=(a:b:c:d).
\]
For each dual point $u$, fix a nonzero representative $(a,b,c,d)\in\Fp^4$ and write
\[
L_u(X,Y,Z,T)=aX+bY+cZ+dT.
\]
Thus $L_u$ defines the plane corresponding to $u$.

Conversely, a point $x=[w]\in\PG(3,p)$ corresponds to the dual plane
\[
x^*=\{u\in\PG^*(3,p):
L_u(w)=0\}.
\]
Under this duality, the planes containing a given line correspond to the points of a line in the dual space, while the planes through a given point correspond to the points of its dual plane.

\subsection{Polynomial evaluation and relations on finite point sets}

Let $D\subseteq\PG(m,K)$ be a finite set, and fix a nonzero representative $a_P\in K^{m+1}$ for each $P\in D$. We denote by $K^D$ the vector space of $K$-valued functions on $D$, and by $K[X_0,\ldots,X_m]_j$ the space of homogeneous polynomials (forms) of degree $j$, with the convention that it is zero for $j<0$.

The evaluation map in degree $j$ is defined by
\[
\operatorname{ev}_{D,j}:
K[X_0,\ldots,X_m]_j\longrightarrow K^D,
\qquad
F\longmapsto (F(a_P))_{P\in D}.
\]
Throughout, we write $F(P)=F(a_P)$ with respect to the chosen representatives.

 The space of evaluation relations is defined by
\begin{equation}\label{eq:relations}
 \mathcal R_j(D)=\left\{c\in K^D:
 \sum_{P\in D}c_PF(a_P)=0\text{ for every }F\in K[X_0,\ldots,X_m]_j\right\}.
\end{equation}

Equivalently, $\mathcal R_j(D)$ is the kernel of the transpose of an evaluation matrix for $\operatorname{ev}_{D,j}$. We say that $D$ is \emph{independent in degree $j$} if $\operatorname{ev}_{D,j}$ is surjective, or equivalently, $\mathcal R_j(D)=0$. The set $D$ is \emph{minimally dependent in degree $j$} if it is dependent in degree $j$ but every proper subset is independent in degree $j$. These notions correspond to $j$-general position and minimal $j$-linkage, respectively, in \cite[Definitions~2.4 and~2.8]{Couvreur}.

If $a_P$ is replaced by $\alpha_Pa_P$, where $\alpha_P\in K^\times$, the corresponding relation coefficients satisfy
\begin{equation}\label{eq:rescale-relation}
c'_P=\alpha_P^{-j}c_P.
\end{equation}

For $c=(c_P)_{P\in D}\in K^D$, we define \(\operatorname{supp}(c) =\{P\in D:c_P\neq0\}.\)  For a form $G$, write $Gc=(G(a_P)c_P)_{P\in D}$. Division of vectors is also coordinatewise and is used only when every denominator is nonzero. For $G\ne0$, we have
\begin{equation}\label{eq:relationmultiply}
G\mathcal R_j(D)\subseteq
\mathcal R_{j-\deg G}(D).
\end{equation}

For a form $G$ and a vector $c\in K^D$, the support of their coordinatewise product is given by
\[
\operatorname{supp}(Gc)
=\{P\in\operatorname{supp}(c):G(a_P)\neq0\}.
\]
We also note that if $c\in\mathcal R_j(D)$ and $\operatorname{supp}(c)$ is independent in degree $j$, then $c=0$.

Moreover, independence in degree $j$ implies independence in every degree at least $j$. Indeed, since $K$ is infinite, we may choose a linear form that does not vanish at any point of $D$, and multiplication by this form induces an injective map $\mathcal R_{j+1}(D)\to\mathcal R_j(D)$.

\subsection{Hilbert functions and plane curves}\label{sec:hilbert}

For finite point sets in the projective plane, the dimensions of the relation spaces can be described using Hilbert functions.

Set $\mathscr S=K[X,Y,Z]$, and let $\mathscr S_v$ denote the space of forms of degree $v$, with $\mathscr S_v=0$ for $v<0$.

{An \emph{ideal} $I\subseteq\mathscr S$ is an additive subgroup closed under multiplication by elements of $\mathscr S$. For $F_1,\ldots,F_r\in\mathscr S$, we write \((F_1,\ldots,F_r) =\left\{\sum_{i=1}^r A_iF_i: A_i\in\mathscr S\right\}\) for the ideal they generate. An ideal is called \emph{homogeneous} if it is generated by homogeneous polynomials. Let $\mathfrak m=(X,Y,Z)$. A homogeneous ideal $I$ is called \emph{saturated} if $\mathfrak m^N A\subseteq I$ for some $N\ge0$ implies $A\in I$, for every $A\in\mathscr S$.

For an ideal $I$ of a polynomial ring $R$ over $K$, the \emph{quotient ring} $R/I$ consists of polynomial classes modulo $I$. Two polynomials represent the same class if their difference belongs to $I$.

Given ideals $I,J\subseteq\mathscr S$, their \emph{ideal quotient} is defined by
\[
I:J=\{A\in\mathscr S:AJ\subseteq I\}.
\]
We write $I:F=I:(F)$ when $F$ is a form. }

For a finite set $D\subseteq\PG(2,K)$, define
\[
(I_D)_v=\{F\in\mathscr S_v:
F(P)=0\text{ for every }P\in D\},
\qquad
I_D=\bigoplus_{v\ge0}(I_D)_v.
\]
Then $I_D$ is a homogeneous ideal.

If $I\subseteq\mathscr S$ is a homogeneous ideal, write $I_v=I\cap\mathscr S_v$. Its \emph{Hilbert function} and \emph{first difference} are defined by
\[
H_I(v)=\dim_K\mathscr S_v-\dim_K I_v,
\qquad
\Delta H_I(v)=H_I(v)-H_I(v-1),
\]
where $H_I(v)=0$ for $v<0$. In particular, write $H_D=H_{I_D}$ and $\Delta H_D=\Delta H_{I_D}$.

By the definition of $I_D$, $\ker(\operatorname{ev}_{D,v})=(I_D)_v$. Hence
\[
\operatorname{rank}(\operatorname{ev}_{D,v})
=\dim_K\mathscr S_v-\dim_K(I_D)_v
=H_D(v).
\]
The transpose of the evaluation matrix has $|D|$ columns, rank $H_D(v)$, and kernel $\mathcal R_v(D)$. Therefore
\[
\dim_K\mathcal R_v(D)=|D|-H_D(v).
\]

\begin{definition}\label{def:zeroset}
For a homogeneous ideal $I\subseteq K[X_0,\ldots,X_m]$, put
\[
V(I)=\{[a]\in\PG(m,K):F(a)=0\text{ for every }F\in I\}.
\]
For forms $F_1,\ldots,F_r$, write $V(F_1,\ldots,F_r)=V((F_1,\ldots,F_r))$. In particular, $V(0)=\PG(m,K)$. When $m=2$ and $F$ is a nonzero form of positive degree, $F=0$ defines a plane curve with zero set $V(F)$; see \cite[Sections~3.1 and~5.1]{Fulton}.
\end{definition}

A nonconstant polynomial in $K[X,Y,Z]$ is \emph{irreducible} if it has no factorization into two nonconstant polynomials. Let
\[
F=c\prod_{i=1}^r F_i^{m_i},
\qquad 0\neq c\in K,\quad m_i\ge1,
\]
where $F_1,\ldots,F_r$ are pairwise coprime irreducible forms over $K$. The \emph{irreducible components} of $V(F)$ are $V(F_1),\ldots,V(F_r)$, with
\[
V(F)=\bigcup_{i=1}^r V(F_i).
\]
The component $V(F_i)$ has multiplicity $m_i$ in the plane curve $F=0$, whose degree is \(\deg F=\sum_{i=1}^r m_i\deg F_i.\)

A plane curve $F=0$ defined over $\Fp$ is \emph{absolutely irreducible} if $F$ is irreducible over $K$.

\subsection{Numerical characters of finite point sets in the plane}

Let $D\subseteq\PG(2,K)$ be a nonempty finite set, and \(s=\min\{v\ge0:(I_D)_v\ne0\}\) be the least degree of a plane curve containing $D$. For $0\le i<s$, define
\[
n_i=1+\max\{v\ge0:\Delta H_D(v)\ge i+1\}.
\]
Thus, $n_i-1$ is the largest degree $v$ for which $\Delta H_D(v)\ge i+1$. These integers satisfy
\[
n_0\ge n_1\ge\cdots\ge n_{s-1}\ge s,
\]
and the sequence $(n_0,\ldots,n_{s-1})$ is called the \emph{numerical character} of $D$ \cite{GP,EP,BallicoNumerical}.

The numerical character determines the first difference of the Hilbert function by
\[
\Delta H_D(v)
=\#\{0\le i<s:i\le v<n_i\}.
\]
In this formula, each index $i$ contributes one precisely for $i\le v\le n_i-1$. Consequently, the numerical character and the Hilbert function determine each other.

The numerical character is called \emph{connected} if $n_{i+1}\ge n_i-1$ for $0\le i<s-1$, and it has a \emph{gap at $j$} if $n_j\le n_{j-1}-2$ for $1\le j<s$. The relevant formulas and decomposition will be recalled in Section \ref{sec:exact}.


\section{Empty planes and projective duality}\label{sec:projection}

Every set satisfying \eqref{eq:hyp} is a cylinder for $p=2,3$ by \cite{BDMS}, and for $p=5,7$ by \cite[Propositions~41 and~52]{KM}, respectively. We therefore assume that $p\ge11$ and that $S$ satisfies \eqref{eq:hyp}.

\subsection{White and black points in the dual projective space}

A projective plane is called \emph{empty} if it is disjoint from $S$. Since $\Pi_\infty\cap S=\varnothing$, hypothesis \eqref{eq:hyp} implies that every projective plane meets $S$ in a multiple of $p$ points. Define the integers $r_\Pi$, the dual point sets $\mathcal{W}$ and $\mathcal B$, and the integer weights $w_{\Pi^*}$ by
\begin{equation}\label{eq:dualdef}
 r_\Pi=\frac{|S\cap\Pi|}{p},\quad
 \mathcal{W}=\{\Pi^*:r_\Pi=0\},\quad
\mathcal B=\{\Pi^*:r_\Pi\ge2\},\quad w_{\Pi^*}=r_\Pi-1\;\;(\Pi^*\in \mathcal B).
\end{equation}
The points of $\mathcal{W}$ and $\mathcal B$ are called \emph{white points} and \emph{black points}, respectively; each black point $b\in\mathcal B$ has positive integer weight $w_b$. The dual points $\Pi^*$ with $r_\Pi=1$ belong to neither set. A similar construction in the planar setting appears in \cite[Section~3]{BMM}.


\subsection{Projection polynomials and bounds for empty planes}

In this subsection, we view $S$ as a subset of $V=\Fp^3$ via the affine coordinates introduced in Section~2.

For a set $S$ that is not a cylinder, we need a polynomial of degree at most $p-2$ representing its intersection sizes with parallel affine lines. The following lemma follows from the classical results recalled in \cite[Lemma~3.2 and Proposition~3.4]{JPRZ}; we use the formulation given in \cite[Propositions~3.1 and~4.1]{KMNS}.

\begin{lemma}[{\cite[Propositions~3.1 and~4.1]{KMNS}}]
\label{lem:KMNS-polynomial}
Let $f:\Fp^3\to\Fp$ be a nonzero function, and $F\in\Fp[X,Y,Z]$ be its unique polynomial representative of degree at most $p-1$ in each variable. Suppose that $\sum_{x\in \Gamma}f(x)=0$ for every affine plane $\Gamma$ in $\Fp^3$. Then $\sum_{x\in\Fp^3}f(x)G(x)=0$ for every polynomial $G\in\Fp[X,Y,Z]$ of total degree at most $p-1$. Moreover, $\deg F\le 2p-3$.
\end{lemma}

\begin{lemma}\label{lem:projection}
Suppose that $S$ satisfies \eqref{eq:hyp} and is not a cylinder. For any direction $L$, choose coordinates such that $\pi_L(x,y,z)=(x,y)$, and define
\[
m(x,y)=|S\cap\pi_L^{-1}(x,y)|,\text{ where } (x,y)\in\Fp^2.
\]
Then there exists a nonzero polynomial $M\in\Fp[X,Y]$ of total degree at most $p-2$ such that
\[
M(x,y)\equiv m(x,y)\pmod p
\]
for all $(x,y)\in\Fp^2$.
\end{lemma}

\begin{proof}
Fix a direction $L$ and the coordinates in the statement. The set $S$ spans $\PG(3,p)$, since otherwise $|S|=p^2$ would make it an affine plane, hence a cylinder. Blokhuis \cite[Lemma~6.2.15]{MattheusThesis}  showed that a $p$-divisible set of $p^2$ points in $\Fp^3$ containing an affine line is a cylinder; see also \cite[Theorem~2]{BDMS}. Since $S$ is not a cylinder, it contains no affine line. Thus $0\le m(x,y)\le p-1$ for every $(x,y)\in\Fp^2$.

Let $\overline{m}(x,y)=m(x,y) \pmod{p}$. Since $\sum_{(x,y)\in\Fp^2}m(x,y)=|S|=p^2$, we have 
\begin{equation}\label{eq:sum_of_m}
\sum_{(x,y)\in\Fp^2}\overline{m}(x,y)=0.
\end{equation}
For every affine line $\ell$ in $\Fp^2$, its inverse image $\pi_L^{-1}(\ell)$ is an affine plane, so
\begin{equation}\label{eq:linesum}
\sum_{(x,y)\in\ell}m(x,y)
=|S\cap\pi_L^{-1}(\ell)|
\equiv0\pmod p.
\end{equation}
Let $M(X,Y)=\sum_{i,j=0}^{p-1}a_{ij}X^iY^j$ be the polynomial representing $\overline{m}(x,y)$, where $a_{ij}\in\Fp$. Since $0\le m(x,y)\le p-1$ and $S$ is nonempty, we have $\overline{m}\neq0$, and hence $M\neq0$.

Define $\widetilde M(X,Y,Z)=M(X,Y)(1-Z^{p-1})$. Then for $(x,y,z)\in\Fp^3$,
\[
\widetilde M(x,y,z)=
\begin{cases}
\overline{m}(x,y), & z=0,\\
0, & z\neq0.
\end{cases}
\]
Every affine plane intersects $Z=0$ in either the whole plane, an affine line, or the empty set; therefore, by \eqref{eq:sum_of_m} and \eqref{eq:linesum}, we have $\sum_{u\in \Gamma}\widetilde M(u)=0$ for every affine plane $\Gamma$ in $\Fp^3$. Since $\widetilde M$ has degree at most $p-1$ in each variable and $M\neq0$, by Lemma \ref{lem:KMNS-polynomial},
\[
\deg M+p-1=\deg\widetilde M\le2p-3.
\]
Thus, $\deg M\le p-2$, as required.
\end{proof}

\begin{lemma}\label{lem:infinity}
Let $S\subseteq\PG(3,p)\setminus\Pi_\infty$ satisfy \eqref{eq:hyp}, and let $\Pi_0$ be any projective plane disjoint from $S$. Then, taking $\Pi_0$ as the plane at infinity, $S$ remains a $p$-divisible set of size $p^2$. Moreover, $S$ is a cylinder with respect to $\Pi_\infty$ if and only if it is a cylinder with respect to $\Pi_0$.
\end{lemma}

\begin{proof}
Let $\Pi_0$ be a projective plane disjoint from $S$. Then the complement of $\Pi_0$ is an affine space containing $S$. Note that every affine plane in this space is of the form $\Pi\setminus\Pi_0$ for some projective plane $\Pi\neq\Pi_0$. Since $S\cap\Pi_0=\varnothing$, we have $|S\cap(\Pi\setminus\Pi_0)|=|S\cap\Pi|$. The first assertion is established.

Now suppose that $S$ is a cylinder with respect to $\Pi_\infty$. Then
\[
S=\bigcup_{i=1}^p(\ell_i\setminus\{v\}),
\]
where $\ell_1,\ldots,\ell_p$ are distinct projective lines through $v\in\Pi_\infty$. We claim that $v\in\Pi_0$. Otherwise, $\Pi_0$ meets some $\ell_i$ in a point other than $v$, which belongs to $S$, contradicting $S\cap\Pi_0=\varnothing$. Hence, all $\ell_i$'s meet $\Pi_0$ at $v$, and so $\ell_i\setminus\{v\}$ ($i=1,\ldots,p$) are parallel affine lines in $\PG(3,p)\setminus\Pi_0$, showing that $S$ is a cylinder with respect to $\Pi_0$. The converse follows by interchanging $\Pi_\infty$ and $\Pi_0$.
\end{proof}

\begin{lemma}\label{lem:planecap}
Suppose that $S$ is not a cylinder. Then for every direction $L$, there are at most $p-2$ affine lines $\ell$ in $\Fp^2\cong V/L$ satisfying $S\cap\pi_L^{-1}(\ell)=\varnothing$ (we call such lines \emph{empty}). Consequently, every dual projective plane $\Sigma$ satisfies
\begin{equation}\label{eq:planecap}
|\mathcal{W}\cap\Sigma|\le p-1.
\end{equation}
\end{lemma}

\begin{proof}
Fix a direction $L$, and let $M\in\Fp[X,Y]$ be the nonzero polynomial given by Lemma~\ref{lem:projection}. Then $\deg M\le p-2$, and $M$ vanishes on every empty line in $V/L$.

Let $\ell$ be an empty line in $V/L$, with equation $g_\ell(X,Y)=0$, where
\[
g_\ell(X,Y)=aX+bY+c
\]
for some $a,b,c\in\Fp$ with $(a,b)\ne(0,0)$. Choose affine coordinates so that $\ell$ is defined by $Y=0$. Since affine changes of coordinates preserve total degree, $M(X,0)$ has degree at most $p-2$. Since $\ell$ is an empty line, we see that $M(x,0)=0$ for every $x\in\Fp$, which implies $M(X,0)$ is the zero polynomial. Hence, $Y$ divides $M$ in these coordinates. Returning to the original coordinates gives $g_\ell|M$. The polynomials $g_\ell$ corresponding to distinct empty lines are pairwise coprime, so their product divides $M$. Since each $g_\ell$ has degree one and $M\neq0$, there are at most $\deg M\le p-2$ empty lines in $V/L$.


Let $\Sigma$ be a plane in the dual projective space. If $\mathcal{W}\cap\Sigma=\varnothing$, the bound is immediate. Otherwise, choose $e_0\in \mathcal{W}\cap\Sigma$. Then $e_0^*$ is an empty plane containing the point $\Sigma^*$. By Lemma~\ref{lem:infinity}, we may take $e_0^*$ as the plane at infinity and identify $\PG(3,p)\setminus e_0^*$ with $V=\Fp^3$, such that the set $S$ is still not a cylinder. Let $L$ be the direction represented by $\Sigma^*$, and let $\pi_L:V\to V/L$ be the quotient map in this affine space.

For each $e\in(\mathcal{W}\cap\Sigma)\setminus\{e_0\}$, put $\Gamma_e=e^*\setminus e_0^*$. Since $e^*\neq e_0^*$, the set $\Gamma_e$ is an affine plane. Write $\Gamma_e=a+U$, where $a\in V$ and $U\le V$ is two-dimensional. The condition $e\in\Sigma$ means that $e^*$ contains $\Sigma^*$. Thus, $\Sigma^*$ lies on the line $e^*\cap e_0^*$ at infinity of $\Gamma_e$, so $L\le U$. Consequently,
\[
\ell_e:=\pi_L(\Gamma_e)=\pi_L(a)+(U/L),
\quad
\pi_L^{-1}(\ell_e)=\Gamma_e.
\]
Since $\dim(U/L)=1$, the set $\ell_e$ is an affine line in $V/L$. Moreover, $e\in \mathcal{W}$ gives $S\cap e^*=\varnothing$, and hence $S\cap\pi_L^{-1}(\ell_e)=S\cap \Gamma_e=S\cap (e^*\setminus e_0^*)=\varnothing$. Therefore, $\ell_e$ is empty. We claim that the map $e \mapsto \ell_e$ is a bijection from $(\mathcal{W}\cap\Sigma)\setminus\{e_0\}$ onto the set of empty affine lines in $V/L$.

If $\ell_e=\ell_f$, then $\Gamma_e=\pi_L^{-1}(\ell_e)=\pi_L^{-1}(\ell_f)=\Gamma_f$. An affine plane determines its containing projective plane uniquely, so $e^*=f^*$ and hence $e=f$. Thus, the map $e\mapsto\ell_e$ is injective.

Conversely, let $\ell$ be an empty affine line in $V/L$. Then $\Gamma=\pi_L^{-1}(\ell)=a+U$ is an affine plane with $L\le U$. Let $\Pi$ be the unique projective plane whose affine part is $\Gamma$. Since the direction of $\Gamma$ contains $L$, the plane $\Pi$ contains $\Sigma^*$. Also, $\Pi\neq e_0^*$. Since $S\cap e_0^*=\varnothing$ and $\ell$ is empty, we have
\[
 S\cap\Pi=S\cap \Gamma
 =S\cap\pi_L^{-1}(\ell)=\varnothing.
\]
It follows that $e=\Pi^*$ belongs to $(\mathcal{W}\cap\Sigma)\setminus\{e_0\}$ and satisfies $\ell_e=\ell$. This proves surjectivity.

Hence the map $e\mapsto\pi_L(e^*\setminus e_0^*)$ is a bijection from $(\mathcal{W}\cap\Sigma)\setminus\{e_0\}$ onto the set of empty affine lines in $V/L$. By the first part of the proof, there are at most $p-2$ such lines. Therefore, $|\mathcal{W}\cap\Sigma|\le 1+(p-2)=p-1$.
\end{proof}

\subsection{Incidence identities and bounds for white points}

In Sections~\ref{sec:projection} and~\ref{sec:exact}, write $n=|\mathcal{W}|$.

\begin{lemma}\label{lem:incidence}
With the notation of \eqref{eq:dualdef}, we have
\begin{equation}\label{eq:energy}
\sum_{b\in\mathcal B}w_b=n-1,\quad
\sum_{b\in\mathcal B}w_b(w_b+1)=p(p-1).
\end{equation}
Moreover, every dual line $\ell$ with $\mathcal{W}\cap\ell\ne\varnothing$ satisfies
\begin{equation}\label{eq:linebalance}
\sum_{b\in\mathcal B\cap\ell}w_b=|\mathcal{W}\cap\ell|-1.
\end{equation}
In particular, $\mathcal B\neq\varnothing$ and $p\le n\le \frac{p^2-p}{2}+1$.
\end{lemma}

\begin{proof}
Let $\Omega$ be the set of planes of $\PG(3,p)$. Then $|\Omega|=p^3+p^2+p+1$. Each point lies in $p^2+p+1$ planes, and each pair of distinct points lies in $p+1$ planes.  By double counting the sizes of $\{ (P,\Pi): \Pi\in \Omega,\, P\in \Pi\cap S \}$ and $\{ (P, Q,\Pi): \Pi\in \Omega,\,\quad P,Q\in \Pi\cap S,\ P\ne Q\}$, we obtain
\begin{align*}
\sum_{\Pi\in \Omega} |S\cap\Pi|&=p^2(p^2+p+1)=p|\Omega|-p,\\
\sum_{\Pi \in \Omega} |S\cap\Pi|(|S\cap\Pi|-1)&=p^2(p^2-1)(p+1).
\end{align*}
By the point-plane duality, an empty plane $\Pi\in\Omega$ corresponds to the white point $\Pi^*\in\mathcal W$. Thus there are exactly $n$ planes with $r_\Pi=0$,  while planes with
$r_\Pi\ge2$ correspond to black points $b=\Pi^*$
with $r_\Pi=w_b+1$. Planes with $r_\Pi=1$ contribute zero to $\sum_\Pi(r_\Pi-1)$ and $\sum_\Pi r_\Pi(r_\Pi-1)$.
Hence, substituting $|S\cap\Pi|=pr_\Pi$ into the preceding identities, we obtain
\begin{align*}
-n+\sum_{b\in\mathcal B}w_b
&=\sum_{\Pi\in\Omega}(r_\Pi-1)=-1,\notag\\
\sum_{b\in\mathcal B}w_b(w_b+1)
&=\sum_{\Pi\in\Omega}r_\Pi(r_\Pi-1)
=p(p-1). \notag
\end{align*}
This proves \eqref{eq:energy}.

Let $L$ be the line of $\PG(3,p)$ dual to $\ell$. 
Since $\mathcal{W}\cap\ell\ne\varnothing$, there exists an empty plane containing $L$,  i.e., $L\cap S=\varnothing$. Each point of $S$ therefore lies in the unique plane through $L$ and there are $p+1$ such planes, and consequently
\[
\sum_{\Pi\supset L}(r_\Pi-1)
=\frac{|S|}{p}-(p+1)=-1.
\]
Among these planes, the empty ones correspond to the white points in $\mathcal W\cap\ell$, and those with $r_\Pi\ge2$ correspond to the black points in $\mathcal B\cap\ell$. Thus $ -|\mathcal W\cap\ell| +\sum_{b\in\mathcal B\cap\ell}w_b=-1,$ which implies \eqref{eq:linebalance}.


The second equality in \eqref{eq:energy} implies $\mathcal B\ne\varnothing$. For every $b\in\mathcal B$, we have
\[
2\le w_b+1=\frac{|S\cap b^*|}{p}\le p.
\]
Summing $2w_b\le w_b(w_b+1)\le pw_b$ over $\mathcal B$ and applying \eqref{eq:energy} gives
\[
2(n-1)\le p(p-1)\le p(n-1).
\]
Therefore $p\le n\le p(p-1)/2+1$. This completes the proof.
\end{proof}

\begin{lemma}\label{lem:cylinder_property}
The following statements hold.
\begin{enumerate}
\item[(i)] Every point of $\PG(3,p)\setminus S$ lies in an empty plane.

\item[(ii)] The white points are coplanar if and only if $S$ is a cylinder.

\item[(iii)] If the black points are collinear, then $S$ is a cylinder.
\end{enumerate}
\end{lemma}

\begin{proof}
\textup{(i)} Let $x\in\PG(3,p)\setminus S$. Each $P\in S$ lies in $p+1$ planes through $x$, so
\[
\sum_{\Pi\ni x}|S\cap\Pi|=(p+1)|S|=p^2(p+1).
\]
If no plane through $x$ is empty, then \eqref{eq:hyp} gives
\[
\sum_{\Pi\ni x}|S\cap\Pi|\ge p(p^2+p+1)>p^2(p+1),
\]
a contradiction. Thus, every point outside $S$ lies in an empty plane, and $S=\PG(3,p) \setminus\bigcup_{e\in \mathcal{W}}e^*$.

\textup{(ii)} Suppose that $\mathcal{W}\subseteq v^*$ for some $v\in\PG(3,p)$. Then $v\in e^*$ for every $e\in \mathcal{W}$. In particular, $v\in\Pi_\infty$, since $\Pi_\infty$ is empty.

For $P\in S$, let $\ell=\langle v,P\rangle$. Since $v\in e^*$ and $P\notin e^*$ for every $e\in \mathcal{W}$, we have $\ell\cap e^*=\{v\}$ for $e\in \mathcal{W}$. The identity $S=\PG(3,p) \setminus\bigcup_{e\in \mathcal{W}}e^*$ in \textup{(i)} therefore gives $\ell\setminus\{v\}\subseteq S$. Since $v\in\Pi_\infty$ and $P$ is arbitrary, $S$ is a disjoint union of parallel affine lines, i.e., $S$ is a cylinder. Conversely, suppose that $S=\bigcup_{i=1}^p(\ell_i\setminus\{v\})$ is a cylinder. A plane not containing $v$ meets each $\ell_i$ at a point of $S$, so every empty plane contains $v$, so $\mathcal{W}\subseteq v^*$.

\textup{(iii)} Suppose that $\mathcal B\subseteq\ell$ for some line $\ell$ in $\PG^*(3,p)$. If $\mathcal{W}\nsubseteq\ell$, choose $e\in \mathcal{W}\setminus\ell$. By Lemma \ref{lem:incidence}, $|\mathcal{W}|\ge p$ so for each $e'\in \mathcal{W}\setminus\{e\}$, \eqref{eq:linebalance} gives
\[
\sum_{b\in\mathcal B\cap\langle e,e'\rangle}
w_b
=|\mathcal{W}\cap\langle e,e'\rangle|-1
\ge1.
\]
Then the line $\langle e,e'\rangle$ intersects $\mathcal{B}$ and hence meets $\ell$ in some point $b=\langle e,e'\rangle\cap\ell$ since $\mathcal{B}\subseteq\ell$. Thus, $e'\in\left\langle e,e' \right\rangle=\left\langle e,b \right\rangle\subseteq\langle e,\ell\rangle$, and so $\mathcal{W}$ lies in the plane $\langle e,\ell\rangle$. If $\mathcal{W}\subseteq\ell$, then $\mathcal{W}$ is also coplanar. In either case, \textup{(ii)} implies that $S$ is a cylinder.
\end{proof}

\begin{lemma}\label{lem:linecap}
If $S$ is not a cylinder, then for every line $\ell$ of $\PG^*(3,p)$ and every $b\in\mathcal B\setminus\ell$, we have
\begin{equation}\label{eq:weightedcap}
|\mathcal{W}\cap\ell|(w_b+1)\le p-1.
\end{equation}
In particular,
\[
|\mathcal{W}\cap\ell|\le\frac{p-1}{2}
\qquad\text{for every line }\ell\subseteq\PG^*(3,p),
\]
and $w_b+1\le(p-1)/2$ for every black point $b\in\mathcal B$.
\end{lemma}

\begin{proof}
Fix a line $\ell$ of $\PG^*(3,p)$ and a point $b\in\mathcal B\setminus\ell$. For each $e\in \mathcal{W}\cap\ell$, \eqref{eq:linebalance} gives $|\mathcal{W}\cap\langle b,e\rangle|\ge w_b+1$. Since $b\notin\ell$, any two distinct points $e,e'\in \mathcal{W}\cap\ell$ satisfy $\langle b,e\rangle\cap\langle b,e'\rangle =\{b\}$. As $b\notin \mathcal{W}$ ($b\in\mathcal{B}$), the sets $\mathcal{W}\cap\langle b,e\rangle$ are pairwise disjoint as $e$ ranges over $\mathcal{W}\cap\ell$. By \eqref{eq:planecap},
\[
|\mathcal{W}\cap\ell|(w_b+1)
\le\sum_{e\in \mathcal{W}\cap\ell}
|\mathcal{W}\cap\langle b,e\rangle|
\le |\mathcal{W}\cap\langle b,\ell\rangle|
\le p-1,
\]
which proves \eqref{eq:weightedcap}.

Since $S$ is not a cylinder, Lemma \ref{lem:cylinder_property} (iii) implies that $\mathcal B$ is not collinear. Therefore, for every line $\ell$, there is a point $b\in\mathcal B\setminus\ell$. As $w_b\ge1$, \eqref{eq:weightedcap} yields $|\mathcal{W}\cap\ell|\le (p-1)/2$.

Finally, for any $b\in\mathcal B$, choose $e\in \mathcal{W}$. Applying \eqref{eq:linebalance} and the preceding bound to $\langle b,e\rangle$ gives $w_b+1\le |\mathcal{W}\cap\langle b,e\rangle| \le (p-1)/2$.
\end{proof}

\subsection{A product identity and an evaluation relation for white points}


Recall that for each $u\in\PG^*(3,p)$, we have fixed a linear form $L_u$. Put $N=\prod_{e\in \mathcal{W}}L_e$, $P=\prod_{b\in\mathcal B}L_b^{w_b}$. Using \eqref{eq:energy} and \eqref{eq:linebalance}, we show that there exist nonzero coefficients $c_e\in\Fp$ $(e\in \mathcal{W})$ such that $\sum_{e\in \mathcal{W}}c_eL_e^{p-2}=0$.

\begin{proposition}\label{prop:product}
There exist coefficients $c_e\in\Fp^{\times}$ $(e\in \mathcal{W})$ such that
\begin{equation}\label{eq:partialfractions}
P=\sum_{e\in \mathcal{W}}c_e\frac{N}{L_e},
\end{equation}
and
\begin{equation}\label{eq:fullrelation}
\sum_{e\in \mathcal{W}}c_eL_e^{p-2}=0.
\end{equation}
With respect to the chosen representatives, \eqref{eq:fullrelation} is equivalent to $(c_e)_{e\in \mathcal{W}}\in\mathcal R_{p-2}(\mathcal{W})$.
\end{proposition}

\begin{proof}
By Lemma \ref{lem:incidence}, $\deg N=n$, $\deg P=\sum_{b\in\mathcal B}w_b=n-1$. Fix $e\in \mathcal{W}$ and consider the restrictions of the polynomials $P$ and $N/L_e$ to the plane $e^*$ defined by $L_e=0$. Since $L_e$ divides neither $P$ nor $N/L_e$, both restrictions are not identically zero.

Choose homogeneous coordinates $(X_0:X_1:X_2)$ on $e^*$ over $\Fp$. For each line $\ell$ through $e$ in $\PG^*(3,p)$, choose a nonzero linear form $\lambda_{e,\ell}\in\Fp[X_0,X_1,X_2]$ defining its dual line $\ell^*$ in $e^*$.

For $u\in\ell\setminus\{e\}$, the incidence relation gives $e^*\cap u^*=\ell^*$. Thus, $L_u|_{e^*}$ is a nonzero scalar multiple of $\lambda_{e,\ell}$. The forms $\lambda_{e,\ell}$ corresponding to distinct lines $\ell$ through $e$ are pairwise coprime.

Recall the definitions $P = \prod_{b\in\mathcal B} L_b^{w_b}$ and $N/L_e = \prod_{u\in \mathcal{W}\setminus\{e\}} L_u$. A linear form $L_b|_{e^*}$ (resp. $L_u|_{e^*}$) divides $\lambda_{e,\ell}$ if and only if $b\in\ell$ (resp. $u\in\ell$). Summing the weights and counting the points of $\mathcal{W}$ on $\ell$, the multiplicities of $\lambda_{e,\ell}$ in $P|_{e^*}$ and $(N/L_e)|_{e^*}$ are therefore $\sum_{b\in\mathcal B\cap\ell}w_b$ and $|\mathcal{W}\cap\ell|-1$, respectively.

By \eqref{eq:linebalance}, these multiplicities are equal for every line $\ell$ through $e$. Since $P$ and $N/L_e$ are both products of linear forms, so are their restrictions on $e^*$. Every factor arises from such a line $\ell$, so unique factorization gives $c_e\in\Fp^{\times}$ such that
\[
P|_{e^*}=c_e (N/L_e)|_{e^*}.
\]
For $e'\neq e$, the polynomial $N/L_{e}$ is divisible by $L_{e'}$. Since 
\[
\left(P-\sum_{e\in \mathcal{W}}c_{e}\frac{N}{L_{e}}\right)|_{e'^*}=P|_{e'^*}-c_{e'}\left(\frac{N}{L_{e'}}\right)|_{e'^*}=0, 
\]
we see that $L_{e'}|\left(P-\sum_{e\in \mathcal{W}}c_{e}\frac{N}{L_{e}}\right)$ for $e'\in \mathcal{W}$. Since the forms $L_{e'}$ are pairwise coprime and $N=\prod_{e\in \mathcal{W}}L_e$, we have 
\[
N | \left(P-\sum_{e\in \mathcal{W}}c_e\frac{N}{L_e}\right).
\]
This polynomial has degree at most $n-1 < \deg N$, so it must be identically zero.

On $K(X,Y,Z,T)$, define the derivation \(\theta =X^p\frac{\partial}{\partial X} +Y^p\frac{\partial}{\partial Y} +Z^p\frac{\partial}{\partial Z} +T^p\frac{\partial}{\partial T},\) where the partial derivatives are formal. For every linear form $L$ over $\Fp$, we have $\theta(L)=L^p$. Put $\Psi =\sum_{b\in\mathcal B}\overline{w}_bL_b^{p-1} -\sum_{e\in \mathcal{W}}L_e^{p-1}$, where $\overline{w}_b=w_b\pmod{p}$. By \eqref{eq:partialfractions}, $\frac{P}{N}=\sum_{e\in \mathcal{W}}c_e\frac{1}{L_e}$. Applying $\theta$ to both sides gives
\[
\theta\left(\frac{P}{N}\right)
=\frac{\theta(P)N-\theta(N)P}{N^2}=\frac{P\Psi}{N}
=-\sum_{e\in \mathcal{W}}c_eL_e^{p-2}.
\]
Since $P\Psi/N$ is a polynomial, $N|P\Psi$. The condition $\mathcal{W}\cap\mathcal B=\varnothing$ gives $\gcd(P,N)=1$, and hence $N|\Psi$. As $\Psi$ has degree at most $p-1 < \deg N$, we obtain $\Psi=0$, and so $\sum_{e\in \mathcal{W}}c_eL_e^{p-2}=0$.

Let $a_e$ be the chosen coefficient vector of $L_e$. For a monomial $F=X^iY^jZ^kT^l$ of degree $p-2$, the coefficient of $F$ in $L_e^{p-2}$ is $(p-2)!F(a_e)/(i!j!k!l!)$. Since $p-2 < p$, $(p-2)!/(i!j!k!l!)$ is nonzero in $\Fp$. Comparing coefficients in \eqref{eq:fullrelation} gives $\sum_{e\in \mathcal{W}}c_eF(a_e)=0$ for every monomial $F$ of degree $p-2$. Since every form $G$ of degree $p-2$ is a linear combination of these monomials, $\sum_{e\in \mathcal{W}}c_eG(a_e)=0$ holds for all $G\in K[X,Y,Z,T]_{p-2}$. Hence, $(c_e)_{e\in \mathcal{W}}\in\mathcal R_{p-2}(\mathcal{W})$.
\end{proof}

\section{Plane intersections and projections}\label{sec:exact}

Let $p\ge11$ be prime, and suppose that $S$ satisfies \eqref{eq:hyp} but is not a cylinder. We prove that every black point has weight one, and hence that all plane intersections with $S$ have size $0$, $p$, or $2p$. For projections of the white point set $\mathcal{W}$ from black points, we then establish the moment relations and the bounds on containing curves used in Section~\ref{sec:ci}.

\subsection{Numerical characters and lower bounds for point sets}

We recall the following formulas for finite sets of plane points in the notation of Section~\ref{sec:hilbert}.

\begin{proposition}[{\cite{EP}, \cite{BallicoNumerical}}]
\label{prop:character-source}
Let $D\subseteq\PG(2,K)$ be a nonempty finite set with numerical character $(n_0,\ldots,n_{s-1})$. Then
\begin{align}
 |D|&=\sum_{i=0}^{s-1}(n_i-i),\label{eq:character-size}\\
 \max\{v\in\mathbb Z:\mathcal R_v(D)\ne0\}&=n_0-2,\label{eq:lastdepend}\\
 \dim_K\mathcal R_v(D)&=\sum_{i=0}^{s-1}
 \bigl((n_i-v-1)_+-(i-v-1)_+\bigr),\label{eq:character-h1}
\end{align}
where $x_+=\max\{x,0\}$.
\end{proposition}

Using $H_D(v)=|D|-\dim_K\mathcal R_v(D)$ and $\dim_K(I_D)_v=\dim_K\mathscr S_v-H_D(v)$, we also obtain
\begin{equation}\label{eq:character-ideal}
 \dim_K(I_D)_v=\dim_K\mathscr S_{v-s}
 +\sum_{i=0}^{s-1}(v+1-n_i)_+
 \qquad(v\ge s-1).
\end{equation}
The formula for $\Delta H_D$ in the definition of the numerical character shows that $H_D$ increases strictly until it reaches $|D|$.

The next proposition gives the gap decomposition for numerical characters. Only its existence assertion is needed here.

\begin{proposition}[{\cite{EP}, \cite[Lemma~8]{BallicoNumerical}}]
\label{prop:gap-source}
Let $D\subseteq\PG(2,K)$ be a nonempty finite set with numerical character $(n_0,\ldots,n_{s-1})$. If $n_j\le n_{j-1}-2$ for some $1\le j<s$, then there is a form $G\in K[X,Y,Z]_j$ such that $D\cap V(G)$ has numerical character $(n_0,\ldots,n_{j-1})$.
\end{proposition}

The next lemma gives the consequence for evaluation relations.

\begin{lemma}\label{lem:gap}
Let $D\subseteq\PG(2,K)$ be a nonempty finite set with numerical character $(n_0,\ldots,n_{s-1})$. Suppose that $n_{j-1}\ge n_j+2$ for some $1\le j<s$. Then there exists a nonzero form $G\in K[X,Y,Z]_j$ such that $D'=D\cap V(G)$ satisfies
\[
|D'|=\sum_{i=0}^{j-1}(n_i-i).
\]
Moreover, suppose that for some $d\ge0$ there exists $c=(c_P)_{P\in D}\in\mathcal R_d(D)$ with $c_P\neq0$ for every $P\in D$.  Then $n_j\ge d+2$.
\end{lemma}

\begin{proof}
By Proposition~\ref{prop:gap-source}, there exists a form $G$ of degree $j$ such that $D'=D\cap V(G)$ has numerical character $(n_0,\ldots,n_{j-1})$. Equation \eqref{eq:character-size} gives \(|D'|=\sum_{i=0}^{j-1}(n_i-i).\) Since $j<s$, the form $G$ cannot vanish on all of $D$, which implies $D'\subsetneq D$.

Let $c=(c_P)_{P\in D}\in\mathcal R_d(D)$ satisfy $c_P\neq0$ for every $P\in D$. Suppose, to the contrary, that $d\ge n_j-1$. Then $d\ge s-1$ and $n_i\le n_j\le d+1$ for $j\le i<s$. Thus \eqref{eq:character-h1} yields
\[
\dim\mathcal R_d(D)
=\sum_{i=0}^{j-1}(n_i-d-1)_+
=\dim\mathcal R_d(D').
\]

For $v=(v_P)_{P\in D'}\in\mathcal R_d(D')$, define $\widetilde v\in K^D$ by
\[
\widetilde v_P=
\begin{cases}
v_P,& P\in D',\\
0,& P\in D\setminus D'.
\end{cases}
\]
For every form $F$ of degree $d$, we have
\[
\sum_{P\in D}\widetilde v_PF(P)
=\sum_{P\in D'}v_PF(P)=0,
\]
so $\widetilde v\in\mathcal R_d(D)$.

The vectors $\widetilde v$ form a subspace of $\mathcal R_d(D)$ of dimension $\dim\mathcal R_d(D')=\dim\mathcal R_d(D)$. Thus there exists $v\in\mathcal R_d(D')$ whose extension by zero is the chosen relation $c$.  Then $c$ satisfies $c_P=0$ for every $P\in D\setminus D'$, contradicting $c_P\neq0$ for every $P\in D$.

\end{proof}

These formulas give the following lower bounds.

\begin{proposition}\label{prop:relationbound}
Let $D\subseteq\PG(2,K)$ be a nonempty finite set, and let $s$ be the least degree of a plane curve containing $D$. Suppose that, for some integer $d\ge0$, there exists $c=(c_P)_{P\in D}\in\mathcal R_d(D)$ with $c_P\neq0$ for every $P\in D$. Then
\begin{equation}\label{eq:truncatedbound}
|D|\ge t(d+3-t)
\qquad(1\le t\le s,\ t\in\mathbb Z).
\end{equation}
\end{proposition}

\begin{proof}
Let $(n_0,\ldots,n_{s-1})$ be the numerical character of $D$. Since $\mathcal R_d(D)\neq0$, \eqref{eq:lastdepend} gives $n_0\ge d+2$. We prove by induction that $n_i\ge d+2-i$ for $0\le i<s$.

Let $1\le i<s$ and suppose that $n_{i-1}\ge d+3-i$. If $n_i\ge n_{i-1}-1$, then \(n_i\ge n_{i-1}-1\ge d+2-i.\) Otherwise, $n_{i-1}\ge n_i+2$, and Lemma~\ref{lem:gap} gives \(n_i\ge d+2\ge d+2-i.\) Thus, since $n_i\ge s$ by definition, we have
\begin{equation}\label{eq:characterlower}
n_i\ge\max\{s,d+2-i\}
\qquad(0\le i<s).
\end{equation}

Let $t$ be an integer with $1\le t\le s$. Since $n_i\ge s>i$ for $0\le i<s$, \eqref{eq:character-size} and \eqref{eq:characterlower} give
\begin{align*}
|D|
=\sum_{i=0}^{s-1}(n_i-i)
\ge\sum_{i=0}^{t-1}(n_i-i)
\ge\sum_{i=0}^{t-1}(d+2-2i)
=t(d+3-t).
\end{align*}
\end{proof}

\subsection{Upper bounds for rational points on plane curves}

For a plane curve $C$ defined over $\Fp$, let $C(\Fp)=C\cap\PG(2,p)$ denote its set of $\Fp$-rational points.

We now turn to upper bounds for points on plane curves, beginning with the following consequence of B\'ezout's theorem.

\begin{proposition}[{\cite[Section~5.3, Corollary~3]{Fulton}}]
\label{prop:bezout-source}
Let $F,G\in K[X,Y,Z]$ be coprime forms of positive degrees $a,b$. Then
\begin{equation}\label{eq:bezout}
 |V(F)\cap V(G)|\le ab.
\end{equation}
\end{proposition}

The following point-count bound extends the result of St\"ohr and Voloch \cite[Theorem~0.1]{SV} to possibly reducible plane curves, as proved by Carlin and Voloch \cite[Theorem~1]{CV}.

\begin{proposition}[{\cite[Theorem~1]{CV}}]
\label{prop:CV}
Let $p$ be prime and let $C$ be a possibly reducible plane curve over $\Fp$ of degree $d<p$. If no irreducible component of $C$ is an $\Fp$-line, then
\[
 |C(\Fp)|\le \frac{d(d+p-1)}2.
\]
\end{proposition}

We shall apply this bound after separating the $\Fp$-line components.

\begin{lemma}\label{lem:polar}
Let $D\subseteq\PG(2,p)$ satisfy $|D\cap\ell|\le(p-1)/2$ for every line $\ell$. If $C$ is a plane curve defined over $\Fp$ of degree $s<p$, then \[
|D\cap C|\le\frac{s(s+p-1)}2.\]
\end{lemma}

\begin{proof}
Let $F\in\Fp[X,Y,Z]_s$ be a form such that $C=V(F)$, and let $\widetilde F$ be the product of the distinct irreducible factors of $F$ over $\Fp$. Write
\[
\widetilde F=L_1\cdots L_rG,
\]
where $L_1,\ldots,L_r$ are its linear factors over $\Fp$ and $G$ has no linear factor over $\Fp$. We have
\[
r+\deg G=\deg\widetilde F\le s<p,
\qquad
C=\bigcup_{i=1}^r V(L_i)\cup V(G).
\]

By hypothesis, \(|D\cap V(L_i)|\le\frac{p-1}{2} (1\le i\le r).\) If $\deg G>0$, Proposition~\ref{prop:CV} gives \(|V(G)(\Fp)|\le\frac{(\deg G)(\deg G+p-1)}2.\) The same inequality holds for $\deg G=0$, since $G$ is a nonzero constant and $V(G)=\varnothing$.

Since $D\subseteq\PG(2,p)$, we have $D\cap V(G)\subseteq V(G)(\Fp)$, and therefore
\begin{align*}
|D\cap C|&\le \sum_{i=1}^r|D\cap V(L_i)|
      +|D\cap V(G)|\\
&\le \frac{r(p-1)}2
      +\frac{(\deg G)(\deg G+p-1)}2
\\
&\le \frac{s(s+p-1)}2.
\end{align*}
\end{proof}

\subsection{Weights of black points and plane intersection sizes}

The next lemma relates independence of a projected set to the sizes of the preimages. A projection line joins the center to a point of the original set; its intersection with the set is the corresponding preimage.

\begin{lemma}\label{lem:bootstrap}
Let $A\subseteq\PG(m,K)$ be a finite set. Fix $b\in\PG(m,K)\setminus A$ and a hyperplane $H\subseteq\PG(m,K)$ with $b\notin H$. The projection from $b$ onto $H$ is
\begin{equation}\label{eq:projection-map}
\pi_{b,H}:A\longrightarrow H,\qquad
P\longmapsto\langle b,P\rangle\cap H.
\end{equation}
Let \(D=\pi_{b,H}(A)=\{P_1,\ldots,P_k\}\) and put \(r_i=|A\cap\langle b,P_i\rangle| (1\le i\le k).\)

Suppose that there exists $c=(c_P)_{P\in A}\in\mathcal R_d(A)$ with $c_P\neq0$ for every $P\in A$. If $D$ is independent in degree $a$, where $0\le a\le d$, then
\[
r_i\ge d-a+2
\qquad(1\le i\le k).
\]
\end{lemma}

\begin{proof}
Choose homogeneous coordinates $(X_0:\cdots:X_{m-1}:T)$ such that
\[
b=(0:\cdots:0:1),
\qquad
H:\ T=0.
\]
For each $P_i\in D$, choose a representative $a_i\in K^m\setminus\{0\}$ such that $P_i=[a_i:0]$. Then
\[
A\cap\langle b,P_i\rangle
=
\{[a_i:t_{i\nu}]:1\le\nu\le r_i\},
\]
where $t_{i1},\ldots,t_{ir_i}$ are pairwise distinct.

Let $c=(c_P)_{P\in A}\in\mathcal R_d(A)$ be the relation in the hypothesis. Let $c_{i\nu}\ne0$ be the corresponding relation coefficients obtained from \eqref{eq:rescale-relation}. Thus \(\sum_{i=1}^k\sum_{\nu=1}^{r_i} c_{i\nu}F(a_i,t_{i\nu})=0\) for every $F\in K[X_0,\ldots,X_{m-1},T]_d$.

For $j\ge0$, define the weighted power sums 
\[
M_i(j)=\sum_{\nu=1}^{r_i}c_{i\nu}t_{i\nu}^{\,j}.
\]
If $0\le j\le d-a$, substituting
$F=T^jG$, with $G\in K[X_0,\ldots,X_{m-1}]_{d-j}$, gives
\begin{equation}\label{eq:bootstrap-evaluation}
\sum_{i=1}^k M_i(j)G(a_i)=0.
\end{equation}
Since $d-j\ge a$ and $D$ is independent in degree $a$, it is independent in degree $d-j$. Hence, for each $i$, there exists such a form $G$ satisfying
\[
G(a_i)=1,
\qquad
G(a_h)=0\quad(h\neq i).
\]
Substituting this choice into \eqref{eq:bootstrap-evaluation} yields
\begin{equation}\label{eq:bootstrap-vanishing}
M_i(j)=0
\qquad
(1\le i\le k,\ 0\le j\le d-a).
\end{equation}

Suppose that $r_i\le d-a+1$ for some $i$. Then \eqref{eq:bootstrap-vanishing} gives
\[
\sum_{\nu=1}^{r_i}
c_{i\nu}t_{i\nu}^{\,j}=0
\qquad
(0\le j<r_i).
\]
The coefficient matrix of this system is \(\bigl(t_{i\nu}^{\,j}\bigr)_{ 0\le j<r_i,\ 1\le\nu\le r_i},\) whose determinant is \(\prod_{1\le\mu<\nu\le r_i} (t_{i\nu}-t_{i\mu})\neq0.\) Hence $c_{i1}=\cdots=c_{ir_i}=0$, contrary to the choice of the relation. Therefore
\[
r_i\ge d-a+2
\qquad(1\le i\le k).
\]
\end{proof}

We now apply these bounds to the white and black point sets $\mathcal{W}$ and $\mathcal B$.

\begin{proposition}\label{prop:simple}
Every black point $b\in\mathcal B$ has weight $w_b=1$. Consequently,
\[
n=\frac{p^2-p+2}{2},\qquad
|S\cap\Pi|\in\{0,p,2p\}
\]
for every projective plane $\Pi$. Moreover, for every $b\in\mathcal B$, there exists a line $\ell$ in $\PG^*(3,p)$ such that
\[
b\in\ell,\qquad |\mathcal{W}\cap\ell|=2.
\]
\end{proposition}

\begin{proof}
Fix $b\in\mathcal B$ and an $\Fp$-rational plane $H\subseteq\PG^*(3,p)$ with $b\notin H$. Using the projection $\pi_{b,H}$ defined in \eqref{eq:projection-map} with $\mathcal{W}$, put \(D=\pi_{b,H}(\mathcal{W}), k=|D|,\) and \(\rho=\min_{P\in D}|\mathcal{W}\cap\langle b,P\rangle|.\) By \eqref{eq:linebalance} and Lemma~\ref{lem:linecap}, we have \(2\le\rho\le\frac{p-1}{2}.\)

Since every projection line contains at least $\rho$ points of $\mathcal{W}$, we have $n\ge\rho k$. For every line $\ell\subseteq H$, \eqref{eq:planecap} gives \(\rho|D\cap\ell| \le |\mathcal{W}\cap\langle b,\ell\rangle| \le p-1.\) Consequently,
\begin{equation}\label{eq:simple-linebound}
k\le\frac{n}{\rho},\qquad
|D\cap\ell|
\le\left\lfloor\frac{p-1}{\rho}\right\rfloor.
\end{equation}

By \eqref{eq:fullrelation}, $\mathcal{W}$ admits a relation in degree $p-2$ with all coefficients nonzero. If $D$ were independent in degree $p-\rho-1$, Lemma~\ref{lem:bootstrap} would imply
\[
|\mathcal{W}\cap\langle b,P\rangle|
\ge\rho+1
\qquad(P\in D),
\]
contrary to the definition of $\rho$. Thus $D$ is dependent in degree $p-\rho-1$.

Choose a minimally dependent subset $D'\subseteq D$ in degree $p-\rho-1$, and let $s$ be the least degree of a curve containing $D'$. Then $\mathcal R_{p-\rho-1}(D')$ contains a relation with all coefficients nonzero. 
Choose $\Fp$-coordinate representatives for the points of $D'$. The evaluation matrix of degree $s$ has entries in $\Fp$, and its rank is unchanged by extending scalars to $K$. Thus there exists a nonzero form $F\in\Fp[X,Y,Z]_s$ vanishing on $D'$. 

If $s=1$, Proposition~\ref{prop:relationbound} and \eqref{eq:simple-linebound} give \(p-\rho+1\le |D'| \le\left\lfloor\frac{p-1}{\rho}\right\rfloor,\) a contradiction.

Suppose that \(2\le s<\frac{p-2\rho+5}{3}.\) Then $s<p$, and Proposition~\ref{prop:relationbound} and Lemma~\ref{lem:polar} yield
\[
|D'|
\ge s(p-\rho+2-s)
>\frac{s(s+p-1)}2
\ge |D'|,
\]
again a contradiction. Therefore
\[
s\ge t:=
\left\lceil\frac{p-2\rho+5}{3}\right\rceil.
\]
Applying \eqref{eq:truncatedbound} to $D'$, we obtain
\begin{equation}\label{eq:minraybound}
n\ge\rho k\ge\rho|D'|
\ge\rho t(p-\rho+2-t).
\end{equation}

Suppose that $\rho\ge3$. Since $t=\lceil(p-2\rho+5)/3\rceil$, we have \(\frac{p-2\rho+5}{3} \le t\le\frac{p-2\rho+7}{3} \le\frac{p-\rho+2}{2}.\) As $z(p-\rho+2-z)$ is increasing for $z\le(p-\rho+2)/2$, \eqref{eq:minraybound} gives
\[
n\ge
\frac{\rho(p-2\rho+5)(2p-\rho+1)}9.
\]
Since $3\le\rho\le(p-1)/2$, we have
\[
\rho(p-2\rho+5)
=3(p-1)+(\rho-3)(p-1-2\rho)
\ge3(p-1).
\]
Together with $2p-\rho+1\ge3(p+1)/2$, this yields
\[
n\ge\frac{p^2-1}{2}>\frac{p^2-p+2}{2},
\]
contrary to $n\le \frac{p^2-p+2}{2}$. Therefore $\rho=2$.

By \eqref{eq:linebalance}, $w_b+1\le |\mathcal{W}\cap\langle b,P\rangle|$ for every $P\in D$. Taking the minimum over $P\in D$ gives $w_b+1\le\rho=2$, and hence $w_b=1$. Since $b$ was arbitrary, we have $w_b=1$ for every $b\in\mathcal B$.

Substituting $w_b=1$ into \eqref{eq:energy} gives $n=\frac{p^2-p+2}{2}$. By the definitions of $\mathcal{W}$ and $\mathcal B$, we have \(|S\cap\Pi|\in\{0,p,2p\}\) for every projective plane $\Pi$.

Finally, since $\rho=2$, there exists $P\in D$ such that $|\mathcal{W}\cap\langle b,P\rangle|=2$. Taking $\ell=\langle b,P\rangle$ completes the proof.
\end{proof}

\subsection{Projection moments and nonvanishing of the first moment}

By Proposition~\ref{prop:simple}, all weights $w_b$ are one. We use the coefficients of \eqref{eq:fullrelation} to define weighted power sums on the lines of a projection.

\begin{proposition}\label{prop:exactmoments}
Let $b\in\mathcal B$, and choose homogeneous coordinates on $\PG^*(3,p)$ such that $b=(0:0:0:1)$. Put $H=\{T=0\}$ and \(D=\pi_{b,H}(\mathcal{W})=\{P_1,\ldots,P_k\}.\) For each $i$, choose $a_i\in\Fp^3\setminus\{0\}$ such that $P_i=[a_i:0]$, and write
\[
\mathcal{W}\cap\langle b,P_i\rangle
=\{[a_i:t_{i\nu}]:1\le\nu\le r_i\},
\]
where $t_{i1},\ldots,t_{ir_i}\in\Fp$ are pairwise distinct.

Let $c_{i\nu}\in\Fp^\times$ be the coefficients supplied by Proposition~\ref{prop:product} for the chosen representatives. For $j\ge0$, put \(M_i(j)=\sum_{\nu=1}^{r_i} c_{i\nu}t_{i\nu}^{\,j}.\) Then
\begin{equation}\label{eq:firstmoment}
M_i(0)=0,\qquad M_i(1)\neq0
\qquad(1\le i\le k).
\end{equation}
In particular, $(M_i(1))_{i=1}^k\in\mathcal R_{p-3}(D)$, with all coordinates nonzero.
\end{proposition}

\begin{proof}
Use the corresponding coordinates $(X:Y:Z:T)$ on $\PG(3,p)$ and normalize the linear forms to the chosen representatives, so that $L_b=T$. For $a_i=(a_{i0},a_{i1},a_{i2})$, put \(\Lambda_i=a_{i0}X+a_{i1}Y+a_{i2}Z.\) The form corresponding to $[a_i:t_{i\nu}]$ is then $\Lambda_i+t_{i\nu}T$.

Recall that \(N=\prod_{e\in \mathcal{W}}L_e, P=\prod_{u\in\mathcal B}L_u^{w_u}\) are the polynomials in \eqref{eq:partialfractions}. Since $L_b=T$ and $w_b=1$ by Proposition~\ref{prop:simple}, put \(P_0=\frac{P}{T} =\prod_{u\in\mathcal B\setminus\{b\}}L_u^{w_u}.\) Then \eqref{eq:partialfractions} gives
\begin{equation}\label{eq:exactexpansion}
\frac{TP_0}{N}
=\sum_{i=1}^k\sum_{\nu=1}^{r_i}
\frac{c_{i\nu}}{\Lambda_i+t_{i\nu}T}.
\end{equation}

Since $N|_{T=0}=\prod_{i=1}^k\Lambda_i^{r_i}\ne0$, setting $T=0$ in \eqref{eq:exactexpansion} gives \(\sum_{i=1}^k\frac{M_i(0)}{\Lambda_i}=0.\) Multiplying by $\prod_{j=1}^k\Lambda_j$ and reducing modulo $\Lambda_i$, we obtain
\[
M_i(0)\prod_{j\ne i}\Lambda_j
\equiv0\pmod{\Lambda_i}.
\]
Since the $\Lambda_i$ are pairwise coprime, it follows that $M_i(0)=0$ for every $i$.

Taking the formal derivative of both sides of \eqref{eq:exactexpansion} with respect to $T$ and setting $T=0$, we obtain
\begin{equation}\label{eq:moment-derivative}
\left.\frac{P_0}{N}\right|_{T=0}
=-\sum_{i=1}^k\frac{M_i(1)}{\Lambda_i^2}.
\end{equation}

By \eqref{eq:linebalance} and $w_b=1$, the multiplicities of $\Lambda_i$ in $N|_{T=0}$ and $P_0|_{T=0}$ are $r_i$ and $r_i-2$, respectively. Multiplying \eqref{eq:moment-derivative} by $N|_{T=0}$ gives
\[P_0|_{T=0}
=-\sum_{j=1}^k M_j(1)\Lambda_j^{r_j-2}
\prod_{h\ne j}\Lambda_h^{r_h}.\]

If $M_i(1)=0$, then $\Lambda_i^{r_i}\mid P_0|_{T=0}$, contrary to the multiplicity $r_i-2$ of $\Lambda_i$ in $P_0|_{T=0}$. Hence $M_i(1)\ne0$ for every $i$.

Finally, applying the evaluation relation \eqref{eq:fullrelation} to $TG$, where $G\in K[X,Y,Z]_{p-3}$, gives
\[
\sum_{i=1}^k M_i(1)G(a_i)=0.
\]
Therefore $(M_i(1))_{i=1}^k\in\mathcal R_{p-3}(D)$.
\end{proof}

For $j\ge0$, put $m_j=(M_i(j))_{i=1}^k\in K^D$. By Proposition~\ref{prop:exactmoments} and the evaluation relation derived from \eqref{eq:fullrelation}, we have
\begin{equation}\label{eq:momentrelations}
m_0=0,\qquad M_i(1)\ne0
\quad(1\le i\le k),
\qquad
m_j\in\mathcal R_{p-2-j}(D)
\quad(0\le j\le p-2).
\end{equation}

We next bound the number of indices $i$ with $r_i\ge3$ and the least degree of a plane curve containing $D$.

\subsection{Curves containing the projected set and their irreducible components}

\begin{corollary}\label{cor:excess}
Let $b\in\mathcal B$, and let $H\subseteq\PG^*(3,p)$ be an $\Fp$-rational plane with $b\notin H$. Put
\[
D=\pi_{b,H}(\mathcal{W})=\{P_1,\ldots,P_k\},
\qquad
r_i=|\mathcal{W}\cap\langle b,P_i\rangle|.
\]
Then
\begin{equation}\label{eq:excess}
\sum_{i=1}^k(r_i-2)=\frac{p^2-p+2}{2}-2k.
\end{equation}
Moreover, $|\{i:r_i\ge3\}|\le \frac{p^2-p+2}{2}-2k$.
\end{corollary}

\begin{proof}
By Proposition~\ref{prop:simple} and \eqref{eq:linebalance}, we have $r_i\ge2$ for every $i$ and $\sum_{i=1}^k r_i=|\mathcal{W}|=\frac{p^2-p+2}{2}$, which gives \eqref{eq:excess}. Each index with $r_i\ge3$ contributes at least one to this sum, so $|\{i:r_i\ge3\}|\le \frac{p^2-p+2}{2}-2k$.
\end{proof}

\begin{proposition}\label{prop:window}
Let $p\ge11$ and $b\in\mathcal B$. Let $H\subseteq\PG^*(3,p)$ be an $\Fp$-rational plane with $b\notin H$, and put $D=\pi_{b,H}(\mathcal{W})$ and $k=|D|$. Let $s$ be the least degree of a plane curve containing $D$. Then
\begin{equation}\label{eq:window}
\left\lceil\frac{p+1}{3}\right\rceil
\le s\le\frac{p-3}{2},
\qquad
s(p-s)\le k\le\frac{p^2-p+2}{4}.
\end{equation}
\end{proposition}

\begin{proof}
By \eqref{eq:momentrelations}, $m_1\in\mathcal R_{p-3}(D)$ has all coefficients nonzero. Since every projection line contains at least two points of $\mathcal{W}$, we have $k\le |\mathcal{W}|/2=(p^2-p+2)/4$.

Suppose that $s\ge(p-1)/2$. Applying Proposition~\ref{prop:relationbound} with $d=p-3$ and $t=(p-1)/2$ gives
\[
k\ge\frac{p^2-1}{4}>\frac{p^2-p+2}{4},
\]
a contradiction. Hence $s\le(p-3)/2<p$.

By \eqref{eq:linebalance} and Proposition~\ref{prop:simple}, every line through $b$ meeting $\mathcal{W}$ contains at least two points of $\mathcal{W}$. Thus, applying \eqref{eq:simple-linebound} with $\rho=2$, we obtain
\begin{equation}\label{eq:projected-linebound}
|D\cap\ell|\le\frac{p-1}{2}
\end{equation}
for every line $\ell\subseteq H$.

Since $D\subseteq\PG(2,p)$, there exists a nonzero form $F\in\Fp[X,Y,Z]_s$ vanishing on $D$. By Proposition~\ref{prop:relationbound} and Lemma~\ref{lem:polar},
\[
s(p-s)\le k\le\frac{s(s+p-1)}2.
\]
Since $s>0$, these inequalities give $2(p-s)\le s+p-1$, and hence $s\ge\lceil(p+1)/3\rceil$. Together with $s\le(p-3)/2$ and $k\le (p^2-p+2)/4$, this proves \eqref{eq:window}.
\end{proof}

We now consider the irreducible components of a curve of least degree containing $D$. The following lemma gives bounds for the points of $D$ lying on each component but on no other component.

\begin{lemma}\label{lem:window-components}
Let $D$ and $s$ be as in Proposition~\ref{prop:window}, and let $F\in\Fp[X,Y,Z]_s$ be a nonzero form vanishing on $D$.

Then $F$ is \emph{squarefree} (that is, $F$ has no repeated irreducible factors over $K$), and every irreducible factor of $F$ over $\Fp$ is absolutely irreducible. Thus, after rescaling, we may write \(F=F_1\cdots F_r,\) where $F_1,\ldots,F_r$ are pairwise coprime absolutely irreducible forms over $\Fp$.

For $1\le i\le r$, put
\[
C_i=V(F_i),\qquad e_i=\deg F_i,\qquad
D_i=(D\cap C_i)\setminus\bigcup_{j\ne i}C_j.
\]
Then $D_i\ne\varnothing$, and $e_i$ is the least degree of a plane curve containing $D_i$. Moreover,
\begin{equation}\label{eq:component-bounds}
e_i(p-s)\le |D_i|
\le\frac{e_i(e_i+p-1)}{2}.
\end{equation}
In particular, $e_i\ge p-2s+1$ for every $i$.
\end{lemma}

\begin{proof}
Suppose that $G^2\mid F$ for some nonconstant irreducible form $G\in K[X,Y,Z]$. Then $F/G$ vanishes on $D$ and has degree $s-\deg G<s$, contradicting the minimality of $s$. Hence $F$ is squarefree.

Suppose that $F_0\mid F$ is irreducible over $\Fp$ but reducible over $K$, and let $Q\mid F_0$ be irreducible over $K$. Define $\operatorname{Fr}:K[X,Y,Z]\to K[X,Y,Z]$ by
\[
\operatorname{Fr}\left(\sum_{i,j,k}a_{ijk}X^iY^jZ^k\right)
=\sum_{i,j,k}a_{ijk}^pX^iY^jZ^k.
\]
Since $F_0$ is irreducible over $\Fp$, there exist $\nu\ge2$ and $c\in K^\times$ such that \(F_0=c\prod_{j=0}^{\nu-1}\operatorname{Fr}^j(Q).\) For every $P\in\PG(2,p)$, we have $\operatorname{Fr}^j(Q)(P)=Q(P)^{p^j}$. Hence $F_0(P)=0$ implies $Q(P)=0$, giving
\[
V(F_0)\cap\PG(2,p)
\subseteq V(Q).
\]

Thus $(F/F_0)Q$ vanishes on $D$, but
\[
\deg((F/F_0)Q)=s-\deg F_0+\deg Q<s,
\]
contrary to the minimality of $s$ over $K$. Hence every irreducible factor of $F$ over $\Fp$ is absolutely irreducible.

We may therefore write
\[
F=F_1\cdots F_r,
\]
where $F_1,\ldots,F_r$ are pairwise coprime absolutely irreducible forms over $\Fp$.

Fix $i$ and put $H_i=F/F_i$. By the definition of $D_i$, we have $H_i(P)=0$ for $P\in D\setminus D_i$ and $H_i(P)\ne0$ for $P\in D_i$. If $D_i=\varnothing$, then $H_i$ vanishes on $D$, contradicting the minimality of $s$, since $\deg H_i=s-e_i<s$.

Suppose that there exists a nonzero form $F_i'\in K[X,Y,Z]$ with $\deg F_i'<e_i$ vanishing on $D_i$. Then $H_iF_i'$ vanishes on $D$ and
\[
\deg(H_iF_i')=s-e_i+\deg F_i'<s,
\]
contrary to the minimality of $s$. Since $D_i\subseteq V(F_i)$, the least degree of a plane curve containing $D_i$ is $e_i$.

By \eqref{eq:momentrelations}, $m_1\in\mathcal R_{p-3}(D)$ has all coefficients nonzero. Since $\deg H_i=s-e_i$, \eqref{eq:relationmultiply} gives
\[
(H_i m_1)|_{D_i}
\in\mathcal R_{p-3-s+e_i}(D_i),
\]
with all coefficients nonzero. Applying Proposition~\ref{prop:relationbound} with $d=p-3-s+e_i$ and $t=e_i$, we obtain $|D_i|\ge e_i(p-s)$.

Since $C_i$ is defined over $\Fp$ and $e_i\le s<p$, \eqref{eq:projected-linebound} and Lemma~\ref{lem:polar} give
\[
|D_i|\le |D\cap C_i|
\le\frac{e_i(e_i+p-1)}2.
\]
Together with $|D_i|\ge e_i(p-s)$, this proves \eqref{eq:component-bounds}. Since $e_i>0$, we obtain $e_i\ge p-2s+1$.
\end{proof}

\section{Interpolation of moments}\label{sec:ci}

We prove that the second and third projection moments are obtained from the first by multiplication by a linear and a quadratic form, respectively.

\begin{theorem}\label{thm:globalmoments}
Let $p\ge11$, and suppose that $S$ satisfies \eqref{eq:hyp} but is not a cylinder. For each $b\in\mathcal B$, let $m_j$ be the moment vectors defined in \eqref{eq:momentrelations}. Then there exist unique forms $L\in\Fp[X,Y,Z]_1$ and $B\in\Fp[X,Y,Z]_2$ such that
\begin{equation}\label{eq:globalmoments}
m_2=Lm_1,\quad m_3=Bm_1.
\end{equation}
\end{theorem}

\subsection{Complete intersections and duality for evaluation relations}

We retain the notation of Section~\ref{sec:hilbert} and first relate evaluation relations to quotients of homogeneous ideals.

Let $F,G\in K[X,Y,Z]$ be coprime forms of positive degrees $d_F,d_G$. Then $V(F,G)$ is finite, and $I=(F,G)$ is a saturated ideal; see \cite[Section~1.4]{EGH}.

For a finite set $D\subseteq V(F,G)$, its \emph{residual ideal} is defined by
\[
I_Z=I:I_D.
\]

This ideal is saturated as well: if $\mathfrak m^N A\subseteq I_Z$, then $\mathfrak m^N AI_D\subseteq I$; the saturation of $I$ gives $AI_D\subseteq I$, and hence $A\in I_Z$.

\begin{proposition}\label{prop:associatedgraded}
Let $F,G\in K[X,Y,Z]$ be coprime forms of positive degrees $d_F,d_G$. Choose a linear form $\ell$ with $V(\ell)\cap V(F,G)=\varnothing$ and coordinates in which $\ell=Z$. Put $f(x,y)=F(x,y,1)$, $g(x,y)=G(x,y,1)$, and $\mathscr T=K[x,y]/(f,g)$.

For $v\ge0$, let $\mathscr T_{\le v}$ be the subspace of $\mathscr T$ consisting of classes represented by polynomials of degree at most $v$, and set $\mathscr T_{\le v}=0$ for $v<0$. Put $f_{d_F}(x,y)=F(x,y,0)$ and $g_{d_G}(x,y)=G(x,y,0)$. Then there is a graded ring isomorphism
\begin{equation}\label{eq:associatedgraded}
\bigoplus_{v\ge0}
\mathscr T_{\le v}/\mathscr T_{\le v-1}
\simeq K[x,y]/(f_{d_F},g_{d_G})
=:\mathscr C.
\end{equation}
On the left, multiplication is defined by
\[
(a+\mathscr T_{\le i-1})
(b+\mathscr T_{\le j-1})
=ab+\mathscr T_{\le i+j-1}
\quad
(a\in\mathscr T_{\le i},\
 b\in\mathscr T_{\le j}).
\]
\end{proposition}

\begin{proof}

For each $v\ge0$, define
\[
\varphi_v:K[x,y]_v\longrightarrow
\mathscr T_{\le v}/\mathscr T_{\le v-1},
\qquad
h\longmapsto\bar h+\mathscr T_{\le v-1},
\]
where $\bar h$ denotes the class of $h$ in $\mathscr T$. Every class in $\mathscr T_{\le v}/\mathscr T_{\le v-1}$ is represented by the homogeneous part of degree $v$ of a polynomial of degree at most $v$, so $\varphi_v$ is surjective.

Since $f=g=0$ in $\mathscr T$ and
\[
\deg(f-f_{d_F})<d_F,\qquad
\deg(g-g_{d_G})<d_G,
\]
we have
\[
f_{d_F}K[x,y]_{v-d_F}
+g_{d_G}K[x,y]_{v-d_G}
\subseteq\ker\varphi_v.
\]
Here and below, $K[x,y]_j=0$ for $j<0$.

Conversely, let $0\ne h\in\ker\varphi_v$. The condition $\varphi_v(h)=0$ means that $\bar h\in\mathscr T_{\le v-1}$, so there exists $q\in K[x,y]$ with $\deg q<v$ and $\bar h=\bar q$. Thus $h-q\in(f,g)$, and we may write
\begin{equation}\label{eq:ideal-representation}
h-q=af+bg
\end{equation}
for some $a,b\in K[x,y]$.

Choose a representation \eqref{eq:ideal-representation} with $\deg q<v$ for which
\[
N=\max\{\deg(af),\deg(bg)\}
\]
is minimal, where $\deg0=-\infty$.

If $N>v$, comparison of the homogeneous parts of degree $N$ in \eqref{eq:ideal-representation} gives
\[
a_{N-d_F}f_{d_F}
+b_{N-d_G}g_{d_G}=0,
\]
where $a_j,b_j$ denote the homogeneous parts of degree $j$ of $a,b$, respectively.

The condition $V(F,G)\cap V(Z)=\varnothing$ implies $\gcd(f_{d_F},g_{d_G})=1$, so there exists a form $c\in K[x,y]_{N-d_F-d_G}$ such that
\[
a_{N-d_F}=cg_{d_G},\qquad
b_{N-d_G}=-cf_{d_F}.
\]
Replacing $a,b$ by $a-cg,b+cf$ preserves \eqref{eq:ideal-representation}, since
\[
(a-cg)f+(b+cf)g=af+bg,
\]
whereas
\[
\deg((a-cg)f)<N,\qquad
\deg((b+cf)g)<N.
\]
This contradicts the minimality of $N$.

It follows that $N\le v$, and comparison in degree $v$ in \eqref{eq:ideal-representation} gives
\[
h=a_{v-d_F}f_{d_F}+b_{v-d_G}g_{d_G}.
\]
Together with the inclusion already proved, this yields
\[
\ker\varphi_v
=f_{d_F}K[x,y]_{v-d_F}
+g_{d_G}K[x,y]_{v-d_G}.
\]

Taking the direct sum of the maps $\varphi_v$ gives a surjective linear map
\[
\begin{aligned}
\varphi:K[x,y]
&\longrightarrow
\bigoplus_{v\ge0}
\mathscr T_{\le v}/\mathscr T_{\le v-1},\\
\sum_{v\ge0}h_v
&\longmapsto
\sum_{v\ge0}
(\overline{h_v}+\mathscr T_{\le v-1}),
\end{aligned}
\]
where $h_v\in K[x,y]_v$. Since $f_{d_F}$ and $g_{d_G}$ are homogeneous, we have
\[
\begin{aligned}
\ker\varphi
&=\bigoplus_{v\ge0}\ker\varphi_v\\
&=\bigoplus_{v\ge0}
\left(
f_{d_F}K[x,y]_{v-d_F}
+g_{d_G}K[x,y]_{v-d_G}
\right)\\
&=(f_{d_F},g_{d_G}).
\end{aligned}
\]

Moreover, for $h\in K[x,y]_i$ and $k\in K[x,y]_j$, we have
\[
\varphi_{i+j}(hk)
=\varphi_i(h)\varphi_j(k).
\]
Thus $\varphi$ is a surjective ring homomorphism with kernel $(f_{d_F},g_{d_G})$, and hence induces the isomorphism
\[
K[x,y]/(f_{d_F},g_{d_G})
\simeq
\bigoplus_{v\ge0}
\mathscr T_{\le v}/\mathscr T_{\le v-1},
\]
which proves \eqref{eq:associatedgraded}.
\end{proof}

For every $v\ge0$, substitution $Z=1$ also induces an isomorphism
\begin{equation}\label{eq:ci-dehom}
\eta_v:\mathscr S_v/I_v\longrightarrow\mathscr T_{\le v},
\qquad Q+I_v\longmapsto Q(x,y,1)+(f,g).
\end{equation}
Surjectivity follows by homogenizing a representative of an element of $\mathscr T_{\le v}$ to degree $v$. To prove injectivity, we first establish $I:Z=I$. If a homogeneous form $Q$ satisfies $ZQ=AF+BG$, we may take $A,B$ homogeneous. Reduction modulo $Z$, together with $\gcd(f_{d_F},g_{d_G})=1$, gives a form $C$ such that $Z\mid A-CG$ and $Z\mid B+CF$. Therefore
\[
ZQ=(A-CG)F+(B+CF)G
\]
can be divided by $Z$ to give $Q\in I$.

Now suppose that $Q(x,y,1)\in(f,g)$. Homogenizing an expression of $Q(x,y,1)$ in the ideal $(f,g)$ gives $Z^NQ\in I$ for some $N\ge0$. Repeated application of $I:Z=I$ then gives $Q\in I$, proving injectivity. Substitution $Z=1$ also preserves products, so the maps $\eta_v$ are compatible with multiplication.

For $v\ge0$, write
\[
\mathscr C_v
=\{h+(f_{d_F},g_{d_G}):
h\in K[x,y]_v\}.
\]
Since $\gcd(f_{d_F},g_{d_G})=1$, we have
\[
f_{d_F}K[x,y]_{v-d_F}
\cap g_{d_G}K[x,y]_{v-d_G}
=f_{d_F}g_{d_G}K[x,y]_{v-d_F-d_G}.
\]
Hence
\[
\begin{aligned}
\dim_K\mathscr C_v
={}&\dim_K K[x,y]_v
-\dim_K K[x,y]_{v-d_F}\\
&-\dim_K K[x,y]_{v-d_G}
+\dim_K K[x,y]_{v-d_F-d_G},
\end{aligned}
\]
where $K[x,y]_j=0$ for $j<0$.

Since $\dim_K K[x,y]_v=v+1$ for $v\ge0$, we have
\[
\sum_{v\ge0}\dim_K K[x,y]_v z^v
=\frac{1}{(1-z)^2},
\]
where $z$ is a formal variable. Therefore,
\[
\begin{aligned}
\sum_{v\ge0}\dim_K\mathscr C_v z^v
&=\frac{(1-z^{d_F})(1-z^{d_G})}
{(1-z)^2}= (\sum_{i=0}^{d_F-1} z^i ) (\sum_{i=0}^{d_G-1} z^i ) \\
\end{aligned}
\]
With $\sigma=d_F+d_G-2$, this gives
\[
\begin{gathered}
\dim_K\mathscr C_\sigma=1,\qquad
\mathscr C_v=0\quad(v>\sigma),\\
\dim_K\mathscr C_i=\dim_K\mathscr C_{\sigma-i}
\quad(0\le i\le\sigma).
\end{gathered}
\]

Since $\mathscr C$ is an Artinian complete intersection, it is Gorenstein by the classical complete-intersection theorem; see \cite[Theorem~CB8]{EGH}. Hence the multiplication pairings
\[
\mathscr C_i\times\mathscr C_{\sigma-i}
\longrightarrow\mathscr C_\sigma,
\qquad (a,b)\longmapsto ab,
\]
are nondegenerate; see \cite[proof of Observation~1 (b)]{DGO}.

The following lemma identifies evaluation relations with homogeneous components of $I_Z/I$. The identifications are compatible with multiplication by forms.

\begin{lemma}\label{lem:ci}
Let $F,G\in\mathscr S$ be coprime forms of positive degrees $d_F,d_G$, and let $D\subseteq V(F,G)$ be finite. Put $I=(F,G)$ and $I_Z=I:I_D$. Choose coordinates as in Proposition~\ref{prop:associatedgraded}. For each $P=(x_P:y_P:1)\in D$, take the representative $a_P=(x_P,y_P,1)$.

Then there are $K$-linear isomorphisms
\begin{equation}\label{eq:cidual}
\Phi_j:\mathcal R_j(D)\longrightarrow
(I_Z/I)_{d_F+d_G-3-j}
\qquad(j\ge0),
\end{equation}
where $(I_Z/I)_v=(I_Z)_v/I_v$. These isomorphisms can be chosen so that
\begin{equation}\label{eq:ci-multiplication}
\Phi_{j-a}(Hc)=H\Phi_j(c)
\end{equation}
for every $c\in\mathcal R_j(D)$ and $H\in\mathscr S_a$, where $0\le a\le j$. Here $Hc=(H(a_P)c_P)_{P\in D}$, and multiplication by $H$ on $I_Z/I$ is induced by multiplication in $\mathscr S$.
\end{lemma}

\begin{proof}
Retain $\mathscr T$ and $\mathscr C$ from Proposition~\ref{prop:associatedgraded}, with $\sigma=d_F+d_G-2$. By \eqref{eq:associatedgraded},
\[
\mathscr T=\mathscr T_{\le\sigma},
\qquad
\dim_K(\mathscr T/\mathscr T_{\le\sigma-1})=1.
\]
We may therefore choose a nonzero linear map $\psi:\mathscr T\to K$ with $\ker\psi=\mathscr T_{\le\sigma-1}$.

We first prove that the pairing $(a,b)\mapsto\psi(ab)$ is nondegenerate. Given $0\ne a\in\mathscr T$, choose $i$ such that $a\in\mathscr T_{\le i}\setminus\mathscr T_{\le i-1}$. Under \eqref{eq:associatedgraded}, the class of $a$ modulo $\mathscr T_{\le i-1}$ corresponds to a nonzero element of $\mathscr C_i$. Nondegeneracy of the multiplication pairing $\mathscr C_i\times\mathscr C_{\sigma-i}\to\mathscr C_\sigma$ therefore gives $b\in\mathscr T_{\le\sigma-i}$ with $ab\notin\mathscr T_{\le\sigma-1}$, and hence $\psi(ab)\ne0$.

More precisely, for $a\in\mathscr T$ and $j\ge0$,
\begin{equation}\label{eq:ci-orthogonal}
a\in\mathscr T_{\le\sigma-j-1}
\quad\Longleftrightarrow\quad
\psi(ab)=0
\text{ for every }b\in\mathscr T_{\le j}.
\end{equation}
The forward implication follows from $\mathscr T_{\le\sigma-j-1}\mathscr T_{\le j} \subseteq\mathscr T_{\le\sigma-1}$. For the converse, suppose that $a\notin\mathscr T_{\le\sigma-j-1}$. For the least $i$ with $a\in\mathscr T_{\le i}$, we then have $i>\sigma-j-1$, or equivalently $\sigma-i\le j$. The nondegeneracy argument above gives $b\in\mathscr T_{\le\sigma-i} \subseteq\mathscr T_{\le j}$ with $\psi(ab)\ne0$.

Since $D\subseteq V(F,G)$, evaluation at each $P\in D$ is well defined on $\mathscr T$. Polynomial interpolation therefore gives a surjection $\mathscr T\to K^D$; denote its kernel by $\mathfrak J$.

For each $c=(c_P)_{P\in D}\in K^D$, nondegeneracy of $(a,b)\mapsto\psi(ab)$ gives a unique $a_c\in\mathscr T$ such that
\begin{equation}\label{eq:ci-pairing}
\psi(a_cb)=\sum_{P\in D}c_Pb(P)
\qquad(b\in\mathscr T).
\end{equation}
For $x\in \mathfrak J$, we have $x(P)=0$ for every $P\in D$, so
\[
\psi((a_cx)b)
=\sum_{P\in D}c_Px(P)b(P)=0
\qquad(b\in\mathscr T).
\]
Nondegeneracy implies $a_cx=0$ for every $x\in \mathfrak J$; thus $a_c\mathfrak J=0$.

Under the substitution $Z=1$, degree-$j$ forms give exactly the evaluations of $\mathscr T_{\le j}$. Consequently, \eqref{eq:relations}, \eqref{eq:ci-pairing}, and \eqref{eq:ci-orthogonal} give
\[
c\in\mathcal R_j(D)
\quad\Longleftrightarrow\quad
a_c\in\mathscr T_{\le\sigma-j-1}.
\]

Fix $j\ge0$ and put $v=\sigma-j-1$. If $v<0$, the preceding equivalence gives $a_c=0$ for every $c\in\mathcal R_j(D)$. Equation~\eqref{eq:ci-pairing} and the surjectivity of evaluation then give $c=0$. Both spaces in \eqref{eq:cidual} are therefore zero, and we take $\Phi_j$ to be the zero map.

Suppose that $v\ge0$, and denote the class of a polynomial $Q$ in the quotient under consideration by $[Q]$. For $c\in\mathcal R_j(D)$, the isomorphism \eqref{eq:ci-dehom} gives a unique class $[Q]\in\mathscr S_v/I_v$ such that $\eta_v([Q])=a_c$. To prove $Q\in I_Z$, note that the image of $I_D$ in $\mathscr T$ under substitution $Z=1$ is exactly $\mathfrak J$: every form in $I_D$ vanishes on $D$, and every polynomial vanishing on $D$ can be homogenized to a form in $I_D$. For any nonzero homogeneous $R\in I_D$, the equality $a_c\mathfrak J=0$ gives
\[
\eta_{v+\deg R}([QR])
=a_c[R(x,y,1)]=0.
\]
Injectivity of $\eta_{v+\deg R}$ gives $QR\in I$. Since $I_D$ is homogeneous, it follows that $QI_D\subseteq I$, so $Q\in I_Z$. We may therefore define
\[
\Phi_j(c)=[Q]\in(I_Z/I)_v.
\]

The map $\Phi_j$ is linear by its construction. Moreover, $\Phi_j(c)=0$ implies $a_c=0$, and then \eqref{eq:ci-pairing} and the surjectivity of $\mathscr T\to K^D$ give $c=0$. Thus $\Phi_j$ is injective.

To prove surjectivity, take $[Q]\in(I_Z/I)_v$ and put $a=\eta_v([Q])$. The inclusion $QI_D\subseteq I$ gives $a\mathfrak J=0$, so the linear functional $b\mapsto\psi(ab)$ vanishes on $\mathfrak J$ and descends to $\mathscr T/\mathfrak J\simeq K^D$. Hence there is a unique $c=(c_P)_{P\in D}\in K^D$ such that
\[
\psi(ab)=\sum_{P\in D}c_Pb(P)
\qquad(b\in\mathscr T).
\]
Uniqueness in \eqref{eq:ci-pairing} gives $a=a_c$. Since $a\in\mathscr T_{\le v}$ and $v=\sigma-j-1$, \eqref{eq:ci-orthogonal} yields $c\in\mathcal R_j(D)$. By construction, $\Phi_j(c)=[Q]$, which proves surjectivity and completes the proof of \eqref{eq:cidual}.

Finally, let $c\in\mathcal R_j(D)$ and $H\in\mathscr S_a$, where $0\le a\le j$. Equation~\eqref{eq:relationmultiply} gives $Hc\in\mathcal R_{j-a}(D)$. Let $h$ be the class of $H(x,y,1)$ in $\mathscr T$. For every $b\in\mathscr T$,
\[
\psi((ha_c)b)
=\psi(a_c(hb))
=\sum_{P\in D}H(a_P)c_Pb(P).
\]
Comparing this equality with \eqref{eq:ci-pairing} and using uniqueness, we obtain $a_{Hc}=ha_c$.

For $v=\sigma-j-1<0$, we have $c=0$, so both sides of \eqref{eq:ci-multiplication} are zero. For $v\ge0$, the compatibility of the maps $\eta_v$ with multiplication gives
\[
\eta_{v+a}\bigl(\Phi_{j-a}(Hc)\bigr)
=a_{Hc}
=ha_c
=\eta_{v+a}\bigl(H\Phi_j(c)\bigr).
\]
The injectivity of $\eta_{v+a}$ proves \eqref{eq:ci-multiplication}.
\end{proof}

\begin{corollary}\label{cor:ciindependence}
Let $F,G\in\mathscr S$ be coprime forms of positive degrees $d_F,d_G$. Then every finite set $Y\subseteq V(F,G)$ is independent in degree $v$ for all $v\ge d_F+d_G-2$.
\end{corollary}

\begin{proof}
If $v\ge d_F+d_G-2$, then $d_F+d_G-3-v<0$. Thus Lemma~\ref{lem:ci}, applied with $D=Y$, gives $\mathcal R_v(Y)=0$, proving the assertion.
\end{proof}

\subsection{Common factors and the first difference of the Hilbert function}

The following lemma describes the dimension increase obtained by multiplying a space of binary forms by linear forms.

\begin{lemma}\label{lem:binary}
Let $\mathscr V\subseteq K[x,y]_d$ be a nonzero subspace with $\dim_K \mathscr V=r$. Then $\dim_K(x\mathscr V+y\mathscr V)\ge r+1$. If equality holds, there exists $0\ne G\in K[x,y]_{d-r+1}$ such that
\[
\mathscr V=GK[x,y]_{r-1}.
\]
\end{lemma}

\begin{proof}
The proof is by induction on $d$. For $r=1$, write $\mathscr V=KG$ with $G\ne0$. Then $x\mathscr V+y\mathscr V=GK[x,y]_1$, which proves both assertions. Suppose that $r\ge2$, so $d\ge1$, and put
\[
\mathscr U=\{H\in K[x,y]_{d-1}:xH,yH\in \mathscr V\}.
\]
The identity $x\mathscr V\cap y\mathscr V=xy\mathscr U$, which follows from $\gcd(x,y)=1$, gives
\begin{equation}\label{eq:binary-dimension}
\dim_K(x\mathscr V+y\mathscr V)=2r-\dim_K \mathscr U.
\end{equation}
If $\mathscr U=0$, then $\dim_K(x\mathscr V+y\mathscr V)=2r>r+1$. If $\mathscr U\ne0$, the induction hypothesis and $x\mathscr U+y\mathscr U\subseteq \mathscr V$ give
\[
\dim_K \mathscr U+1\le\dim_K(x\mathscr U+y\mathscr U)\le r,
\]
and hence $\dim_K(x\mathscr V+y\mathscr V)\ge r+1$ by \eqref{eq:binary-dimension}.

If $\dim_K(x\mathscr V+y\mathscr V)=r+1$, then \eqref{eq:binary-dimension} gives $\dim_K \mathscr U=r-1$. The preceding inequalities now give $\dim_K(x\mathscr U+y\mathscr U)=r=\dim_K \mathscr V$, so $\mathscr V=x\mathscr U+y\mathscr U$. Applying the equality case of the induction hypothesis to $\mathscr U$, we obtain
\[
\mathscr U=GK[x,y]_{r-2}
\qquad\text{for some }0\ne G\in K[x,y]_{d-r+1}.
\]
Therefore, $ \mathscr V=x\mathscr U+y\mathscr U =G\bigl(xK[x,y]_{r-2}+yK[x,y]_{r-2}\bigr) =GK[x,y]_{r-1}.$
\end{proof}

The next lemma is a form of Davis's common-factor theorem \cite{Davis85}; see also \cite[Theorem~2.3]{GHM}. We include a proof using Lemma~\ref{lem:binary}.

\begin{lemma}\label{lem:plateaufactor}
Let $I\subseteq\mathscr S$ be a saturated homogeneous ideal with finite projective zero set, and let $z$ be the least degree of a nonzero form in $I$. Suppose that
\[
q\ge z-1,\qquad
\Delta H_I(q)=\Delta H_I(q+1)=c>0.
\]
Then a form of degree $c$ divides every form in $I$ of degree at most $q+1$.
\end{lemma}

\begin{proof}
If $q=z-1$, then $I_{z-1}=0$ and $\Delta H_I(z-1)=\Delta H_I(z)=z=c$. Consequently,
\[
\dim_K I_z=z+1-\Delta H_I(z)=1.
\]
A nonzero generator of $I_z$ has degree $z=c$ and divides every form in $I$ of degree at most $q+1=z$.

Suppose that $q\ge z$. For $v=q,q+1$, let $G_v$ be the greatest common divisor of the nonzero forms in $I_v$, up to a nonzero scalar. After division by $G_v$, these forms have only finitely many common projective zeros. Choose a line $\ell=0$ avoiding all these points and $V(I)$, and which is not a component of either $G_q$ or $G_{q+1}$. Let $\mathscr V_v$ be the space of restrictions of $I_v$ to this line. Since $I$ is saturated and $\ell$ does not vanish on $V(I)$, we have $I:\ell=I$; hence
\[
\dim_K \mathscr V_v
=\dim_K I_v-\dim_K I_{v-1}
=v+1-\Delta H_I(v).
\]
In particular,
\[
\dim_K \mathscr V_q=q+1-c>0,\qquad
\dim_K \mathscr V_{q+1}=q+2-c.
\]
Using homogeneous coordinates $x,y$ on the line, we have $x\mathscr V_q+y\mathscr V_q\subseteq \mathscr V_{q+1}$. Lemma~\ref{lem:binary} therefore gives
\[
q+2-c\le\dim_K(x\mathscr V_q+y\mathscr V_q)\le\dim_K \mathscr V_{q+1}=q+2-c.
\]
Its equality case shows that $\mathscr V_q$ and $\mathscr V_{q+1}=x\mathscr V_q+y\mathscr V_q$ have the same greatest common divisor, up to a nonzero scalar, of degree $q-\dim_K \mathscr V_q+1=c$.

By the choice of $\ell$, the restrictions of the forms obtained after division by $G_v$ have no common zero on the line. Their greatest common divisor is therefore one, so the greatest common divisor of $\mathscr V_v$ is the restriction of $G_v$. Since $\ell$ is not a component of $G_v$, this restriction has degree $\deg G_v$, and thus $\deg G_q=\deg G_{q+1}=c$. Moreover, $\mathscr S_1 I_q\subseteq I_{q+1}$ implies that $G_{q+1}$ divides $XA,YA,ZA$ for every $A\in I_q$. As $\gcd(X,Y,Z)=1$, it follows that $G_{q+1}\mid A$ and hence $G_{q+1}\mid G_q$. The two forms therefore differ by a nonzero scalar.

It remains to consider $A\in I_v$ with $v<q$. For every monomial $M$ of degree $q+1-v$, we have $MA\in I_{q+1}$, so $G_q\mid MA$. These monomials have greatest common divisor one; hence $G_q\mid A$. Thus $G_q$ divides every form in $I$ of degree at most $q+1$.
\end{proof}

\subsection{Moment recurrences and rational functions}

Retain the projection notation and representatives of Proposition~\ref{prop:exactmoments}:
\[
D=\pi_{b,H}(\mathcal{W})=\{P_1,\ldots,P_k\},\qquad
b=(0:0:0:1),\qquad H=\{T=0\}.
\]
For each $i$, the chosen representative $a_i\in\Fp^3\setminus\{0\}$ gives
\[
\mathcal{W}\cap\langle b,P_i\rangle
=\{[a_i:t_{i\nu}]:1\le\nu\le r_i\},
\qquad
M_i(j)=\sum_{\nu=1}^{r_i}c_{i\nu}t_{i\nu}^{\,j},
\]
where the parameters $t_{i\nu}\in\Fp$ are pairwise distinct and $c_{i\nu}\in\Fp^\times$ are the coefficients of the evaluation relation.

Fix $i$ with $r_i=2$, and write $t_{i1}=\alpha$ and $t_{i2}=\beta$. Equation~\eqref{eq:firstmoment} gives $c_{i1}=-c_{i2}\ne0$, and therefore $M_i(j)=c_{i1}(\alpha^j-\beta^j)$. For $j\ge0$, define
\begin{equation}\label{eq:twopoint}
h_j(i)=\frac{M_i(j+1)}{M_i(1)}
=\frac{\alpha^{j+1}-\beta^{j+1}}{\alpha-\beta}
=\sum_{\nu=0}^j\alpha^{j-\nu}\beta^\nu.
\end{equation}
In particular,
\[
h_0(i)=1,\qquad h_1(i)=\alpha+\beta,\qquad
h_2(i)=\alpha^2+\alpha\beta+\beta^2,\qquad
h_1(i)^2-h_2(i)=\alpha\beta.
\]
Both $\alpha$ and $\beta$ satisfy $T^2-(\alpha+\beta)T+\alpha\beta=0$. Substitution in \eqref{eq:twopoint} therefore gives
\begin{equation}\label{eq:tworecurrence}
h_j(i)=h_1(i)h_{j-1}(i)
-\bigl(h_1(i)^2-h_2(i)\bigr)h_{j-2}(i)
\qquad(j\ge2).
\end{equation}

The next lemma gives conditions under which rational functions satisfying this recurrence are polynomials.

\begin{lemma}\label{lem:poles}
Let $0\ne A\in\mathscr S$ have degree $d_A\ge1$, and let $h_0,\ldots,h_{2d_A+2}$ be rational functions satisfying $Ah_n\in\mathscr S$ for $0\le n\le2d_A+2$. Suppose that $h_0=1$, $h_1=f$, and $h_2=g$, and put $\gamma =f^2-g$. If
\begin{equation}\label{eq:polerecurrence}
h_n=fh_{n-1}-\gamma h_{n-2}
\end{equation}
holds for $2\le n\le d_A+2$, then $f\in\mathscr S$. If it holds for $2\le n\le2d_A+2$, then $f,g\in\mathscr S$.
\end{lemma}

\begin{proof}
Suppose first that \eqref{eq:polerecurrence} holds for $2\le n\le d_A+2$. The initial identity $h_2h_0-h_1^2=-\gamma $ and the calculation
\[
h_{n+1}h_{n-1}-h_n^2
=\gamma \bigl(h_nh_{n-2}-h_{n-1}^2\bigr)
\qquad(2\le n\le d_A+1)
\]
give, by induction,
\begin{equation}\label{eq:cassini}
h_{n+1}h_{n-1}-h_n^2=-\gamma ^n
\qquad(1\le n\le d_A+1).
\end{equation}
Since $Ah_j\in\mathscr S$ for all indices involved, multiplication by $A^2$ yields
\begin{equation}\label{eq:cassini-polynomial}
A^2\gamma ^n=(Ah_n)^2-(Ah_{n+1})(Ah_{n-1})\in\mathscr S
\qquad(1\le n\le d_A+1).
\end{equation}

Fix an irreducible factor $\pi$ of $A$. For $0\ne P\in\mathscr S$, let $v_\pi(P)$ be the multiplicity of $\pi$ in $P$, and extend this notation by
\[
v_\pi(P/Q)=v_\pi(P)-v_\pi(Q),\qquad v_\pi(0)=+\infty.
\]
For rational functions $u,v$, we have
\[
v_\pi(uv)=v_\pi(u)+v_\pi(v),\qquad
v_\pi(u+v)=\min\{v_\pi(u),v_\pi(v)\}
\quad\text{if }v_\pi(u)\ne v_\pi(v).
\]
The assumptions $Ah_n\in\mathscr S$ and $v_\pi(A)\le d_A$ give
\begin{equation}\label{eq:poles-bound}
v_\pi(h_n)\ge-v_\pi(A)\ge-d_A
\qquad(0\le n\le2d_A+2).
\end{equation}
Taking $n=d_A+1$ in \eqref{eq:cassini-polynomial}, we obtain
\[
2v_\pi(A)+(d_A+1)v_\pi(\gamma )\ge0,\qquad
v_\pi(\gamma )\ge-\frac{2v_\pi(A)}{d_A+1}
\ge-\frac{2d_A}{d_A+1}>-2.
\]
As $v_\pi(\gamma )\in\mathbb Z\cup\{+\infty\}$, it follows that $v_\pi(\gamma )\ge-1$.

If $v_\pi(f)<0$, then $v_\pi(\gamma )>2v_\pi(f)$. Induction using \eqref{eq:polerecurrence} gives $v_\pi(h_n)=nv_\pi(f)$ for $0\le n\le d_A+1$: the cases $n=0,1$ follow from $h_0=1$ and $h_1=f$, while
\[
v_\pi(fh_{n-1})
=nv_\pi(f)
<v_\pi(\gamma )+(n-2)v_\pi(f)
=v_\pi(\gamma h_{n-2})
\]
proves the induction step. In particular,
\[
v_\pi(h_{d_A+1})=(d_A+1)v_\pi(f)\le-d_A-1,
\]
contradicting \eqref{eq:poles-bound}. Thus $v_\pi(f)\ge0$ for every irreducible factor $\pi$ of $A$. Because $Af\in\mathscr S$, every factor of a reduced denominator of $f$ would divide $A$; hence unique factorization gives $f\in\mathscr S$.

Suppose now that \eqref{eq:polerecurrence} holds for $2\le n\le2d_A+2$. The same induction extends \eqref{eq:cassini} and \eqref{eq:cassini-polynomial} to $1\le n\le2d_A+1$. Taking $n=2d_A+1$, we obtain
\[
v_\pi(\gamma )\ge-\frac{2v_\pi(A)}{2d_A+1}
\ge-\frac{2d_A}{2d_A+1}>-1,
\]
so $v_\pi(\gamma )\ge0$ for every irreducible factor $\pi$ of $A$. Since $f\in\mathscr S$ and $Ag\in\mathscr S$,
\[
A\gamma =Af^2-Ag\in\mathscr S.
\]
Unique factorization therefore gives $\gamma \in\mathscr S$, and $g=f^2-\gamma \in\mathscr S$.
\end{proof}

\subsection{Interpolation of moments on an irreducible plane curve}\label{sec:momentplateau}

The preceding results allow the recurrence on a sufficiently large subset of an irreducible curve to determine the first two moment identities on the whole set.

\begin{lemma}\label{lem:momentplateau}
Let $D_0\subseteq\PG(2,K)$ be a nonempty finite set, and let $e$ be the least degree of a plane curve containing $D_0$. Suppose that $D_0$ lies on an irreducible plane curve $C$ of degree $e$. Let $d_0,m$ be integers satisfying $d_0\ge2e$, $m\ge1$, and $4m+2<e$. Suppose that the numerical character $(n_0,\ldots,n_{e-1})$ of $D_0$ is connected. Put $b_i=n_i-(d_0+2-i)$ for $0\le i<e$, and assume that
\begin{equation}\label{eq:plateauhyp}
b_{3m}=b_{3m+1}=b_{3m+2}=b_{3m+3}=m.
\end{equation}

Let
\[
\mu_n\in\mathcal R_{d_0-n}(D_0)
\qquad(0\le n\le2m+2),
\]
with $(\mu_0)_P\ne0$ for every $P\in D_0$, and define
\[
h_n(P)=\frac{(\mu_n)_P}{(\mu_0)_P}
\qquad(P\in D_0,\ 0\le n\le2m+2).
\]
Suppose that $D'\subseteq D_0$ satisfies $|D'|>e(4m+2)$ and
\[
h_n(P)=h_1(P)h_{n-1}(P)
-\bigl(h_1(P)^2-h_2(P)\bigr)h_{n-2}(P)
\]
for every $P\in D'$ and $2\le n\le2m+2$. Then there exist forms $L\in\mathscr S_1$ and $B\in\mathscr S_2$ such that
\[
\mu_1=L\mu_0,\qquad \mu_2=B\mu_0.
\]
\end{lemma}

\begin{proof}
Connectedness of the numerical character gives $n_i-1\le n_{i+1}\le n_i$, and hence
\[
0\le b_{i+1}-b_i=n_{i+1}-n_i+1\le1
\qquad(0\le i<e-1).
\]
Put $R=b_{e-1}$. Its definition gives
\[
n_{e-1}=d_0+3-e+R,\qquad e+n_{e-1}-3=d_0+R.
\]
Using $b_{3m+3}=m$ from \eqref{eq:plateauhyp}, we obtain
\begin{equation}\label{eq:Rrange}
m=b_{3m+3}\le R
\le b_{3m+3}+(e-3m-4)=e-2m-4.
\end{equation}

Choose an irreducible form $F$ of degree $e$ with $C=V(F)$. For $v=n_{e-1}$, the last summand in \eqref{eq:character-ideal} equals one, so
\[
\dim_K(I_{D_0})_{n_{e-1}}
>\dim_K\mathscr S_{n_{e-1}-e}
=\dim_K(F\mathscr S_{n_{e-1}-e}).
\]
As $F\mathscr S_{n_{e-1}-e}\subseteq(I_{D_0})_{n_{e-1}}$, there exists
\[
G\in(I_{D_0})_{n_{e-1}}\setminus F\mathscr S_{n_{e-1}-e}.
\]
Thus $\deg G=n_{e-1}$ and $\gcd(F,G)=1$. Put $I=(F,G)$ and $I_Z=I:I_{D_0}$, and choose coordinates and representatives as in Lemma~\ref{lem:ci}. By \eqref{eq:rescale-relation}, replacing $a_P$ by $\alpha_Pa_P$ changes $h_n(P)$ to $\alpha_P^nh_n(P)$. All terms of the recurrence have the same degree in $\alpha_P$, and the identities $\mu_1=L\mu_0$, $\mu_2=B\mu_0$ transform in the same way. It therefore suffices to use these representatives.

For $0\le n\le2m+2$, equation~\eqref{eq:Rrange} gives
\[
R+n\le e-2<\min\{e,n_{e-1}\},\qquad I_{R+n}=0,
\]
where $n_{e-1}=d_0+3-e+R>e$. Since $e+n_{e-1}-3=d_0+R$, the isomorphism \eqref{eq:cidual} identifies $\mathcal R_{d_0-n}(D_0)$ with $(I_Z)_{R+n}$. Define
\[
Q_n:=\Phi_{d_0-n}(\mu_n)\in(I_Z)_{R+n}
\qquad(0\le n\le2m+2).
\]
In particular, $Q_0=\Phi_{d_0}(\mu_0)\ne0$. As $R+n<d_0$, the multiplication rule gives
\[
Q_0\mu_n,\ Q_n\mu_0\in\mathcal R_{d_0-R-n}(D_0).
\]
By \eqref{eq:ci-multiplication}, both have image $[Q_0Q_n]$ under the isomorphism
\[
\mathcal R_{d_0-R-n}(D_0)\simeq(I_Z/I)_{2R+n}.
\]
Its injectivity therefore yields
\begin{equation}\label{eq:Qcross}
Q_0\mu_n=Q_n\mu_0
\qquad(0\le n\le2m+2).
\end{equation}

The integers $i-b_i$ are nondecreasing, and \eqref{eq:plateauhyp} gives
\[
(3m+j)-b_{3m+j}=2m+j\qquad(0\le j\le3).
\]
For $v=2m+1,2m+2$, the indices satisfying $i-b_i\le v$ are therefore exactly $0,\ldots,v+m$. Consequently, \eqref{eq:character-h1} and \eqref{eq:cidual} give
\[
\begin{aligned}
\dim_K(I_Z)_{R+v}-\dim_K(I_Z)_{R+v-1}
&=\sum_{i=0}^{e-1}
 \bigl((v+1-i+b_i)_+-(v-i+b_i)_+\bigr)\\
&=\#\{0\le i<e:i-b_i\le v\}
=v+m+1.
\end{aligned}
\]
Here $I_{R+v}=I_{R+v-1}=0$ and $d_0-v\ge e-1$ by \eqref{eq:Rrange} and $d_0\ge2e$. It follows that
\begin{equation}\label{eq:plateaudim}
\begin{aligned}
\Delta H_{I_Z}(R+v)
&=R+v+1-\#\{0\le i<e:i-b_i\le v\}\\
&=R-m
\qquad(v=2m+1,2m+2).
\end{aligned}
\end{equation}
The ideal $I_Z$ is saturated, and $V(I_Z)\subseteq V(F,G)$ is finite. If $R>m$, then $0\ne Q_0\in(I_Z)_R$ and \eqref{eq:plateaudim} allow an application of Lemma~\ref{lem:plateaufactor} with $q=R+2m+1$. The resulting form $G_0$ has degree $R-m$ and divides every $Q_n$, since $Q_n\in(I_Z)_{R+n}$ and $R+n\le q+1$. If $R=m$, take $G_0=1$. In either case,
\[
Q_n=G_0N_n,\qquad
N_n\in\mathscr S_{m+n},\qquad N_0\ne0
\quad(0\le n\le2m+2).
\]

For $0\le n\le2m+2$, put $w_n=N_0\mu_n-N_n\mu_0$. Equations \eqref{eq:relationmultiply} and \eqref{eq:Qcross} give
\[
w_n\in\mathcal R_{d_0-n-m}(D_0),\qquad
G_0w_n=Q_0\mu_n-Q_n\mu_0=0.
\]
Thus $(w_n)_P=0$ whenever $G_0(P)\ne0$, and
\[
\operatorname{supp}(w_n)
\subseteq D_0\cap V(G_0)\subseteq V(F,G_0).
\]
If $R=m$, then $G_0=1$ and $w_n=0$. If $R>m$, irreducibility of $F$ and $\deg G_0=R-m<e$ give $\gcd(F,G_0)=1$. Moreover,
\[
d_0-n-m\ge e+(R-m)-2,
\]
because $d_0\ge2e$ and $R+n\le e-2$. Corollary~\ref{cor:ciindependence} makes $\operatorname{supp}(w_n)$ independent in this degree. The restriction of $w_n$ to its support is therefore zero. Hence, in both cases,
\begin{equation}\label{eq:Across}
N_0\mu_n=N_n\mu_0
\qquad(0\le n\le2m+2).
\end{equation}
Since every coordinate of $\mu_0$ is nonzero, this implies
\[
N_n(P)=N_0(P)h_n(P)
\qquad(P\in D_0,\ 0\le n\le2m+2).
\]
In particular, $N_0(P)=0$ implies $N_n(P)=0$ for every $n$.

If a nonzero form of degree at most $4m+2$ vanished on $D'$, its degree would be less than $e$, so it would be coprime to $F$. Equation~\eqref{eq:bezout} would then give $|D'|\le e(4m+2)$, contrary to the hypothesis. Thus every form of degree at most $4m+2$ vanishing on $D'$ is zero. Put $L=N_1/N_0$ and $B=N_2/N_0$. For $2\le n\le m+2$, we obtain
\begin{equation}\label{eq:recurrenceclear1}
N_0^2N_n-N_0N_1N_{n-1}
+(N_1^2-N_0N_2)N_{n-2}=0.
\end{equation}
Indeed, the left-hand side is a form of degree $3m+n\le4m+2$. At $P\in D'$, substitution of $N_j(P)=N_0(P)h_j(P)$ gives
\[
N_0(P)^3\bigl(
h_n(P)-h_1(P)h_{n-1}(P)
+(h_1(P)^2-h_2(P))h_{n-2}(P)
\bigr)=0.
\]
The preceding degree bound therefore proves \eqref{eq:recurrenceclear1}. Division by $N_0^3$ yields
\[
\frac{N_n}{N_0}
=L\frac{N_{n-1}}{N_0}
-(L^2-B)\frac{N_{n-2}}{N_0}
\qquad(2\le n\le m+2).
\]
The rational functions $N_n/N_0$ have initial terms $1,L,B$, and multiplication by $N_0\in\mathscr S_m$ makes each of them a polynomial. Lemma~\ref{lem:poles} therefore gives $L\in\mathscr S$. As $N_0L=N_1\in\mathscr S_{m+1}$, homogeneity gives $L\in\mathscr S_1$.

With $L$ now a polynomial, for $2\le n\le2m+2$ we have
\begin{equation}\label{eq:recurrenceclear2}
N_0N_n-N_0LN_{n-1}
+(N_0L^2-N_2)N_{n-2}=0.
\end{equation}
The left-hand side has degree $2m+n\le4m+2$. If $P\in D'$ and $N_0(P)\ne0$, then $L(P)=h_1(P)$ and $N_j(P)=N_0(P)h_j(P)$, so its value is
\[
N_0(P)^2\bigl(
h_n(P)-h_1(P)h_{n-1}(P)
+(h_1(P)^2-h_2(P))h_{n-2}(P)
\bigr)=0.
\]
If $N_0(P)=0$, then all $N_j(P)=0$, and its value is again zero. The same consequence of \eqref{eq:bezout} therefore proves \eqref{eq:recurrenceclear2}. Division by $N_0^2$ gives \eqref{eq:polerecurrence} for $N_n/N_0$ through index $2m+2$. The second assertion of Lemma~\ref{lem:poles} yields $B\in\mathscr S$; since $N_0B=N_2\in\mathscr S_{m+2}$, we obtain $B\in\mathscr S_2$.

Finally, \eqref{eq:Across}, together with $N_0L=N_1$ and $N_0B=N_2$, gives
\[
\begin{gathered}
\mu_1-L\mu_0\in\mathcal R_{d_0-1}(D_0),\qquad
\mu_2-B\mu_0\in\mathcal R_{d_0-2}(D_0),\\
N_0(\mu_1-L\mu_0)=N_0(\mu_2-B\mu_0)=0.
\end{gathered}
\]
Both relations are supported on $D_0\cap V(N_0)\subseteq V(F,N_0)$. Since $\deg N_0=m<e$, the forms $F,N_0$ are coprime, and
\[
d_0-1\ge d_0-2\ge2e-2\ge e+m-2.
\]
Corollary~\ref{cor:ciindependence} therefore gives $\mu_1=L\mu_0$ and $\mu_2=B\mu_0$.
\end{proof}

\subsection{Interpolation of moments on each irreducible component}\label{sec:components}

The coordinates and representatives are those of Proposition~\ref{prop:exactmoments}. After applying Lemma~\ref{lem:momentplateau}, these representatives are restored. The moment identities are homogeneous and are therefore preserved by this change.

Let $F=F_1\cdots F_r$, with $C_i,D_i$ and $e_i$ as in Lemma~\ref{lem:window-components}, and put
\[
J=D\setminus\bigsqcup_{i=1}^r D_i,
\qquad \delta=p-2s\ge3,
\qquad \lambda_i=|D_i|-e_i(p-s).
\]
Thus $J$ consists of the points of $D$ lying on at least two components. Proposition~\ref{prop:window} and Lemma~\ref{lem:window-components} give
\begin{align}
e_i&\ge\delta+1\ge4,\label{eq:componentdegree}\\
\lambda_i&\ge0,\label{eq:lambdapos}\\
k-s(p-s)=|J|+\sum_i\lambda_i
&\le\frac{\delta^2-p+2}{4},\label{eq:budget}\\
\lambda_i&\le\frac{e_i(e_i-\delta-1)}2.\label{eq:componentbudget}
\end{align}
Indeed, $k=|J|+\sum_i|D_i|$ and $\sum_i e_i=s$ give $k-s(p-s)=|J|+\sum_i\lambda_i$. Since $k\le (p^2-p+2)/4$,
\[
k-s(p-s)\le\frac{p^2-p+2}{4}-s(p-s)
=\frac{p^2-p+2-4sp+4s^2}{4}
=\frac{(p-2s)^2-p+2}{4}
=\frac{\delta^2-p+2}{4}.
\]
The upper bound in \eqref{eq:component-bounds} also gives
\[
\lambda_i
\le\frac{e_i(e_i+p-1)}2-e_i(p-s)
=\frac{e_i(e_i-\delta-1)}2,
\]
which proves \eqref{eq:componentbudget}.

The following numerical lemma will be applied to the bounds \eqref{eq:budget} and \eqref{eq:componentbudget}.

\begin{lemma}\label{lem:numeric}
Let $e\ge4$ be an integer and $z\ge2e+3$ a real number. Put
\[
a=\left\lfloor\frac{e+1}{4}\right\rfloor,
\qquad \tau(e)=ae-\frac32a(a+1).
\]
Then $
e(z-e)+\tau(e)>
\min\left\{\frac{z^2-z+2}{4},\frac{e(e+z-1)}2\right\}.$
\end{lemma}

\begin{proof}
Since $(e-2)/4\le a\le(e+1)/4$ and $x\mapsto ex-3x(x+1)/2$ is concave,
\[
\tau(e)\ge\min\left\{
\frac{5e^2-16e+12}{32},
\frac{5e^2-10e-15}{32}\right\}
\ge\frac{e(5e-16)}{32}.
\]
Suppose, to the contrary, that
\[
e(z-e)+\tau(e)\le\frac{z^2-z+2}{4},
\qquad
e(z-e)+\tau(e)\le\frac{e(e+z-1)}2.
\]
Put $v=z-2e\ge3$. The second inequality gives
\[
v\le e-1-\frac{2\tau(e)}e\le\frac{11e}{16},
\]
whereas the first inequality gives
\[
v^2-v\ge4\tau(e)+2e-2\ge\frac{5e^2-16}{8}.
\]
As $x\mapsto x^2-x$ is increasing for $x\ge3$, it follows that
\[
\frac{5e^2-16}{8}
\le v^2-v
\le\frac{121e^2}{256}-\frac{11e}{16}.
\]
This is impossible, since
\[
\frac{5e^2-16}{8}
-\left(\frac{121e^2}{256}-\frac{11e}{16}\right)
=\frac{39e^2}{256}+\frac{11e}{16}-2>0
\qquad(e\ge4).
\]
\end{proof}

The moment identities on each $D_i$ follow from Lemma~\ref{lem:momentplateau} and these numerical bounds.

\begin{proposition}\label{prop:componentmoments}
For every $i$, there exist forms $L_i\in\mathscr S_1$ and $B_i\in\mathscr S_2$ such that
\[
m_2=L_im_1,\qquad m_3=B_im_1
\quad\text{on }D_i.
\]
\end{proposition}

\begin{proof}
Fix $i$, and write $e=e_i$ and $H_i=F/F_i$. Put
\[
d_0=p-s+e-3,
\qquad \mu_n=(H_i m_{n+1})|_{D_i}
\quad(0\le n\le p-3),
\]
and
\[
(D_i)_2=\{P\in D_i:|\mathcal{W}\cap\langle b,P\rangle|=2\}.
\]
The relevant sets in Lemma~\ref{lem:momentplateau} are $D_0=D_i$, $C=C_i$, and $D'=(D_i)_2$.

The form $H_i$ is nonzero on $D_i$ and vanishes on $D\setminus D_i$. For $0\le n\le p-3$, equation~\eqref{eq:momentrelations} gives $m_{n+1}\in\mathcal R_{p-3-n}(D)$. Since $\deg H_i=s-e$, the multiplication rule \eqref{eq:relationmultiply}, followed by restriction to $D_i$, gives
\[
\mu_n\in\mathcal R_{p-3-n-(s-e)}(D_i)
=\mathcal R_{d_0-n}(D_i)
\]
whenever the resulting degree is nonnegative. Moreover,
\[
(\mu_0)_P=H_i(P)(m_1)_P\ne0
\quad(P\in D_i),
\qquad
\frac{\mu_n}{\mu_0}=\frac{m_{n+1}}{m_1}
\quad\text{on }D_i.
\]
Thus these ratios satisfy \eqref{eq:tworecurrence} on $(D_i)_2$. The absolute irreducibility of $C_i$ and the minimality of its degree $e$ follow from Lemma~\ref{lem:window-components}.

The degree condition in Lemma~\ref{lem:momentplateau} follows from
\[
d_0-2e=(s-e)+\delta-3\ge0.
\]
If $m\ge1$ and $4m+2<e$, then
\[
2m+2<e\le s\le\frac{p-3}{2}.
\]
Hence all $0\le n\le2m+2$ satisfy $n+1\le p-2$ and $d_0-n\ge2e-n>0$. In particular, all moment relations required by Lemma~\ref{lem:momentplateau} are available.

Let $(n_0,\ldots,n_{e-1})$ be the numerical character of $D_i$. Suppose that it is not connected, and let $j<e$ be its last gap. Since $\mu_0\in\mathcal R_{d_0}(D_i)$ has no zero coordinate, Lemma~\ref{lem:gap} gives
\[
n_j\ge d_0+2,
\qquad n_\ell\ge n_{j-1}\ge n_j+2\ge d_0+4
\quad(0\le\ell<j).
\]
The degree-$j$ curve supplied by that lemma therefore contains at least
\[
\sum_{\ell=0}^{j-1}(n_\ell-\ell)
\ge j(d_0+4)-\frac{j(j-1)}2>je
\]
points of $D_i$, where the strict inequality follows from $d_0\ge2e$ and $j\le e-1$. Its defining form is coprime to $F_i$, because $F_i$ is irreducible and $j<e$. Equation~\eqref{eq:bezout} then bounds the number of common points by $je$, a contradiction. Thus the numerical character of $D_i$ is connected.

For $0\le j<e$, put $b_j=n_j-(d_0+2-j)$. Equation~\eqref{eq:characterlower}, applied to $\mu_0$ in degree $d_0$, gives $b_j\ge0$. Connectedness and monotonicity give $n_j-1\le n_{j+1}\le n_j$, and hence $0\le b_{j+1}-b_j\le1$. Consequently,
\[
\begin{gathered}
b_j\ge0\quad(0\le j<e),
\qquad
b_j\le b_{j+1}\le b_j+1\quad(0\le j<e-1),\\
\sum_{j=0}^{e-1}b_j=\lambda_i.\notag
\end{gathered}
\]
The last identity follows from \eqref{eq:character-size}, since
\[
|D_i|
=\sum_{j=0}^{e-1}(n_j-j)
=e(d_0+3-e)+\sum_{j=0}^{e-1}b_j
=e(p-s)+\sum_{j=0}^{e-1}b_j.
\]

Suppose, for a contradiction, that no pair $L_i,B_i$ satisfies both asserted moment identities. If $b_3=0$, then $b_0=b_1=b_2=b_3=0$. The bound $b_{j+1}\le b_j+1$ also gives $b_j\le j-3$ for $j\ge3$. For $v=0,1,2$, the second positive-part term in \eqref{eq:character-h1} vanishes because $j-(d_0-v)-1<0$ for $0\le j<e$. Therefore
\[
\dim_K\mathcal R_{d_0-v}(D_i)
=\sum_{j=0}^{e-1}(v+1-j+b_j)_+
=\sum_{j=0}^{v}(v+1-j)
=\binom{v+2}{2}
\qquad(v=0,1,2).
\]
Here the terms with $j\ge3$ vanish because $j-b_j\ge3$, and those with $v<j\le2$ vanish because $b_j=0$.

For $v=1,2$, multiplication by $\mu_0$ defines a linear map
\[
\mathscr S_v\longrightarrow\mathcal R_{d_0-v}(D_i),
\qquad A\longmapsto A\mu_0.
\]
This map is injective: $A\mu_0=0$ implies $A|_{D_i}=0$, whereas $\deg A=v<e$ and $e$ is the least degree of a curve containing $D_i$. Both spaces have dimension $\binom{v+2}{2}$, so the map is an isomorphism. Since $\mu_1\in\mathcal R_{d_0-1}(D_i)$ and $\mu_2\in\mathcal R_{d_0-2}(D_i)$, there exist $L_i\in\mathscr S_1$ and $B_i\in\mathscr S_2$ with
\[
\mu_1=L_i\mu_0,\qquad \mu_2=B_i\mu_0.
\]
Cancelling $H_i(P)\ne0$ at every $P\in D_i$ gives $m_2=L_im_1$ and $m_3=B_im_1$ on $D_i$, contrary to the assumption. Thus $b_3\ge1$.

Suppose that $m\ge1$, $4m+2<e$, and $b_{3j}\ge j$ for $1\le j\le m$. Since $3m+3\le4m+2<e$, the entry $b_{3m+3}$ is defined. If $b_{3m+3}\le m$, monotonicity gives
\[
b_{3m}=b_{3m+1}=b_{3m+2}=b_{3m+3}=m.
\]
Under this equality, the remaining cardinality condition in Lemma~\ref{lem:momentplateau} follows from the stronger inequality $|(D_i)_2|>e(5m+2)$, proved below.

By \eqref{eq:excess}, at most $\frac{p^2-p+2}{2}-2k$ points $P\in D$ satisfy $|\mathcal{W}\cap\langle b,P\rangle|\ge3$. Since every projection line contains at least two points of $\mathcal{W}$, $|(D_i)_2|\ge|D_i|-(\frac{p^2-p+2}{2}-2k)$. In addition, $k=|J|+\sum_j|D_j|$ and \eqref{eq:lambdapos} imply $k\ge|D_i|+(s-e)(p-s)$. Hence
\begin{equation}\label{eq:tworaycomponent}
|(D_i)_2|\ge|D_i|+2k-\frac{p^2-p+2}{2}
\ge3|D_i|+2(s-e)(p-s)-\frac{p^2-p+2}{2}.
\end{equation}
For $1\le j\le m$, the inequalities $b_{3j}\ge j$ and $b_\ell\ge b_{3j}$ for $\ell\ge3j$ give
\[
\lambda_i=\sum_{\ell=0}^{e-1}b_\ell
\ge\sum_{j=1}^{m}\bigl|\{\ell:b_\ell\ge j\}\bigr|
\ge\sum_{j=1}^{m}(e-3j)
=me-\frac32m(m+1).
\]
Substituting $|D_i|=e(p-s)+\lambda_i$ into \eqref{eq:tworaycomponent} therefore gives
\begin{align*}
|(D_i)_2|-e(5m+2)
&\ge(2s+e)(p-s)+3\lambda_i-e(5m+2)-\frac{p^2-p+2}{2}\\
&\ge(2s+e)(p-s)-2me-\frac92m(m+1)-2e-\frac{p^2-p+2}{2}.
\end{align*}

The expression $2me+9m(m+1)/2$ is increasing for $m\ge0$. Since $4m+2<e$ implies $m\le(e-3)/4$,
\[
2me+\frac92m(m+1)
\le\frac{25e^2-66e-27}{32}.
\]
Combining these inequalities gives
\[
|(D_i)_2|-e(5m+2)
\ge\frac{32(2s+e)(p-s)-25e^2+2e-16p^2+16p-5}{32}.
\]
Substitution of $p=2s+\delta$ rewrites the numerator as
\begin{equation}\label{eq:positivenumerator}
7e^2+32(s-e)(e+1)+16\delta(2e-\delta)+34e+16\delta-5.
\end{equation}
By \eqref{eq:componentdegree}, $2e-\delta\ge\delta+2>0$, and $s-e\ge0$. All terms in \eqref{eq:positivenumerator} other than $-5$ are therefore nonnegative, while $7e^2+34e+16\delta-5>0$. Thus
\[
|(D_i)_2|>e(5m+2)>e(4m+2).
\]
All hypotheses of Lemma~\ref{lem:momentplateau} are now satisfied. It gives both moment identities after cancellation of $H_i$ on $D_i$, contrary to the assumption. Therefore
\[
b_{3j}\ge j\ (1\le j\le m),\quad 4m+2<e
\quad\Longrightarrow\quad b_{3m+3}\ge m+1.
\]

Put $a=\lfloor(e+1)/4\rfloor$. For $1\le m\le a-1$,
\[
4m+2\le4a-2\le e-1,
\qquad 3m+3\le3a\le e-1.
\]
For the last inequality, if $e\ge7$ then $3a\le3(e+1)/4\le e-1$. If $e=4,5,6$, then $a=1$ and $3a=3\le e-1$, and the range $1\le m\le a-1$ is empty. Starting from $b_3\ge1$, the preceding implication therefore gives $b_{3j}\ge j$ for every $1\le j\le a$. The same counting argument yields
\begin{equation}\label{eq:surpluslower}
\lambda_i\ge\sum_{j=1}^a(e-3j)=\tau(e).
\end{equation}

Set $z=\delta+2e$. Since $\delta\ge3$ and $e\le s$,
\[
2e+3\le z\le2s+\delta=p.
\]
Equation~\eqref{eq:surpluslower} gives $e(z-e)+\tau(e)\le e(z-e)+\lambda_i$. The nonnegativity in \eqref{eq:lambdapos} and the bounds \eqref{eq:budget}--\eqref{eq:componentbudget} give, respectively,
\begin{align*}
e(z-e)+\lambda_i&\le e(e+\delta)+\frac{\delta^2-p+2}{4}
=\frac{z^2-p+2}{4}\le\frac{z^2-z+2}{4},\\
e(z-e)+\lambda_i&\le e(e+\delta)+\frac{e(e-\delta-1)}2
=\frac{e(e+z-1)}2.
\end{align*}
Together, these inequalities give
\[
e(z-e)+\tau(e)\le
\min\left\{\frac{z^2-z+2}{4},\frac{e(e+z-1)}2\right\},
\]
contradicting Lemma~\ref{lem:numeric}. Hence the required forms $L_i,B_i$ exist.
\end{proof}

\subsection{Moment identities on the whole projected set}

Proposition~\ref{prop:componentmoments} gives forms $L_i,B_i$ satisfying the moment identities on each $D_i$. It remains to prove that these forms are independent of $i$ and that the identities hold on $J$. The latter assertion will follow from the next independence result.

\begin{lemma}\label{lem:intersectionindependence}
Let $F_1,\ldots,F_r$ be pairwise coprime homogeneous plane forms of positive degrees $e_i$, with $\sum_i e_i=s$. Every finite set of distinct points lying on at least two of the curves $V(F_i)$ is independent in every degree $v\ge s-2$.
\end{lemma}

\begin{proof}
The proof is by induction on $r$. For $r=1$, the set is empty. Suppose that $r\ge2$, let $J$ be a finite set of distinct points lying on at least two curves, and put
\[
J_0=J\cap V(F_r),\qquad J_1=J\setminus J_0,
\qquad G=F_1\cdots F_{r-1}.
\]
Then $J_0\subseteq V(F_r,G)$ and $\gcd(F_r,G)=1$. Since $\deg F_r+\deg G=s$, Corollary~\ref{cor:ciindependence} gives independence of $J_0$ in degree $s-2$. Every point of $J_1$ lies on at least two of the first $r-1$ curves, so the induction hypothesis gives independence of $J_1$ in degree $s-e_r-2$.

For arbitrary values $c_P\in K$, $P\in J$, choose $A\in\mathscr S_{s-2}$ satisfying $A(P)=c_P$ on $J_0$. Since $F_r(P)\ne0$ for $P\in J_1$, independence of $J_1$ gives $B\in\mathscr S_{s-e_r-2}$ such that
\[
B(P)=\frac{c_P-A(P)}{F_r(P)}
\qquad(P\in J_1).
\]
If $J_1=\varnothing$, take $B=0$. The form $A+F_rB\in\mathscr S_{s-2}$ satisfies
\[
(A+F_rB)(P)=
\begin{cases}
A(P)=c_P,&P\in J_0,\\
A(P)+F_r(P)B(P)=c_P,&P\in J_1.
\end{cases}
\]
Thus $J$ is independent in degree $s-2$. The property of independence recalled in Section~\ref{sec:prelim} then gives independence in every degree $v\ge s-2$.
\end{proof}

\begin{proof}[Proof of Theorem~\ref{thm:globalmoments}]
By Proposition~\ref{prop:componentmoments}, choose forms $L_i\in\mathscr S_1$ and $B_i\in\mathscr S_2$ satisfying $m_2=L_im_1$ and $m_3=B_im_1$ on $D_i$. Put $A_{i,1}=L_i$ and $A_{i,2}=B_i$, and fix $\nu\in\{1,2\}$.

Suppose that $A_{i,\nu}\ne A_{j,\nu}$ for some $i\ne j$, and define
\[
H_{ij}=F/(F_iF_j),
\qquad w=H_{ij}(m_{\nu+1}-A_{j,\nu}m_1).
\]
Equations \eqref{eq:momentrelations} and \eqref{eq:relationmultiply} give
\[
m_{\nu+1}-A_{j,\nu}m_1\in\mathcal R_{p-3-\nu}(D),
\qquad \deg H_{ij}=s-e_i-e_j,
\]
and hence
\[
w\in\mathcal R_{p-3-\nu-s+e_i+e_j}(D).
\]
The identity $m_{\nu+1}=A_{j,\nu}m_1$ on $D_j$ gives $w|_{D_j}=0$, while $H_{ij}$ vanishes on every $C_\ell$ with $\ell\notin\{i,j\}$. Therefore
\[
Y:=\operatorname{supp}(w)
\subseteq D_i\cup(D\cap C_i\cap C_j)
\subseteq C_i\cap(D_i\cup J).
\]
This inclusion also accounts for intersections with a third component, where $H_{ij}=0$.

For $P\in D_i$, the moment identity on $D_i$ gives
\[
w(P)=H_{ij}(P)(A_{i,\nu}-A_{j,\nu})(P)(m_1)_P.
\]
Both $H_{ij}(P)$ and $(m_1)_P$ are nonzero. Moreover, $A_{i,\nu}-A_{j,\nu}$ is a nonzero form of degree $\nu<e_i$, and is therefore coprime to the irreducible form $F_i$. Equation~\eqref{eq:bezout} and the lower bound in \eqref{eq:component-bounds} give
\[
|Y|\ge|D_i|-\nu e_i
\ge e_i(p-s-\nu)>e_i(e_i-1).
\]
The strict inequality follows from $p-s=s+\delta\ge e_i+3$ and $\nu\le2$. The defining form of any curve of degree less than $e_i$ is coprime to $F_i$, so \eqref{eq:bezout} bounds its number of common points with $C_i$ by $e_i(e_i-1)$. Since $Y\subseteq C_i$, the least degree of a curve containing $Y$ is exactly $e_i$.

By the definition of $Y$, the restriction $w|_Y$ is a relation in degree $p-3-\nu-s+e_i+e_j$ with every coefficient nonzero. Proposition~\ref{prop:relationbound}, applied with $t=e_i$, gives
\[
|Y|\ge e_i\bigl((p-3-\nu-s+e_i+e_j)+3-e_i\bigr)
=e_i(p-s+e_j-\nu).
\]
On the other hand,
\[
|Y|\le|D_i|+|J|=e_i(p-s)+\lambda_i+|J|.
\]
Combining these inequalities with \eqref{eq:lambdapos} and \eqref{eq:componentdegree} yields
\[
k-s(p-s)\ge\lambda_i+|J|
\ge e_i(e_j-\nu)
\ge(\delta+1)(\delta-1)=\delta^2-1.
\]
This contradicts \eqref{eq:budget}, because
\[
(\delta^2-1)-\frac{\delta^2-p+2}{4}
=\frac{3\delta^2+p-6}{4}>0.
\]
Thus $A_{i,\nu}=A_{j,\nu}$ for all $i,j$ and $\nu\in\{1,2\}$. Denote the common linear and quadratic forms by $L$ and $B$.

The moment identities hold on every $D_i$, so the relations
\[
m_2-Lm_1\in\mathcal R_{p-4}(D),
\qquad m_3-Bm_1\in\mathcal R_{p-5}(D)
\]
are supported on $J$. Every point of $J$ lies on at least two components, and \eqref{eq:window} gives $p-4\ge p-5\ge s-2$. Lemma~\ref{lem:intersectionindependence} therefore makes both relations zero. If $r=1$, then $J=\varnothing$ and the same conclusion is immediate. Hence both identities hold on all of $D$.

For uniqueness, the bound $s\ge\lceil(p+1)/3\rceil\ge4$ gives $(I_D)_1=(I_D)_2=0$. Since $(m_1)_P\ne0$ for every $P\in D$, two choices of $L$ differ by an element of $(I_D)_1$, and two choices of $B$ differ by an element of $(I_D)_2$. Thus $L$ and $B$ are unique.

Finally, the representatives $a_i$ and the vectors $m_1,m_2,m_3$ have coordinates in $\Fp$. The equations $Lm_1=m_2$ and $Bm_1=m_3$ are therefore linear systems over $\Fp$ in the coefficients of $L$ and $B$. Row reduction over $\Fp$ shows that their unique solutions over $K$ have coefficients in $\Fp$. This proves the theorem.
\end{proof}


\section{Proof of the strong cylinder conjecture}\label{sec:quadric}

Let $p\ge11$, and suppose that $S$ satisfies \eqref{eq:hyp} but is not a cylinder. By Proposition~\ref{prop:simple}, we have $|\mathcal{W}|=1+p(p-1)/2$, and every black point $b\in\mathcal B$ has weight $w_b=1$. We first show that no quadric contains $\mathcal{W}$, using the incidence bounds and the classification of generalized hyperfocused arcs. Here a quadric is the zero set of a nonzero homogeneous quadratic form over $K$; it may be reducible or singular.

We then use Theorem~\ref{thm:globalmoments} to construct a quadric containing all points of $\mathcal{W}$ on projection lines that meet $\mathcal{W}$ in exactly two points. An independence result for the remaining projected points shows that this quadric contains all of $\mathcal{W}$. This contradiction completes the proof.

We recall the projection notation of Proposition~\ref{prop:exactmoments}. Whenever a black point $b\in\mathcal B$ is chosen as the projection center, we use the coordinates and representatives from that proposition. Thus $b=(0:0:0:1)$, $H=\{T=0\}$, and $D=\pi_{b,H}(\mathcal{W})=\{P_1,\ldots,P_k\}$. For each $i$, we have $P_i=[a_i:0]$ and
\[
\mathcal{W}\cap\langle b,P_i\rangle
=\{[a_i:t_{i\nu}]:1\le\nu\le r_i\},
\]
where $a_i\in\Fp^3\setminus\{0\}$, $r_i=|\mathcal{W}\cap\langle b,P_i\rangle|$, and $t_{i1},\ldots,t_{ir_i}\in\Fp$ are pairwise distinct. We also retain the coefficients $c_{i\nu}$ from Proposition~\ref{prop:exactmoments} and the moment vectors $m_j$ from \eqref{eq:momentrelations}.

\subsection{Excluding quadrics containing the white point set}

A \emph{$k$-arc} in $\PG(2,p)$ is a set of $k$ points with no three collinear. A line containing two points of the arc is called a \emph{secant}. A \emph{blocking set} for a family of lines is a point set meeting every line in the family. Following Giulietti and Montanucci \cite[Sections~1--2]{GM}, an arc $A\subseteq\PG(2,p)$ is called \emph{generalized hyperfocused} if there exists a set $\mathcal T\subseteq\PG(2,p)\setminus A$ of size $|A|-1$ that meets every secant of $A$, where the points of $\mathcal T$ need not be collinear; see also \cite[Section~1]{BMM}.

\begin{theorem}[Blokhuis--Marino--Mazzocca {\cite[Theorem~7]{BMM}}]\label{thm:hyperfocused}
Let $p$ be an odd prime. A generalized hyperfocused arc in $\PG(2,p)$ with at least two points has size two or four.
\end{theorem}

We first apply Theorem~\ref{thm:hyperfocused} to show that $\mathcal{W}$ contains three collinear points.

\Needspace{8\baselineskip}
\begin{lemma}\label{lem:nocap}

The set $\mathcal{W}$ contains three collinear points.

\end{lemma}

\begin{proof}
Suppose that no three white points of $\mathcal{W}$ are collinear. Let $\Sigma$ be any dual plane containing at least two points of $\mathcal{W}$, and put $A=\mathcal{W}\cap\Sigma$. Then $A$ is an arc. Since $w_b=1$ for every $b\in\mathcal B$, equation~\eqref{eq:linebalance} shows that every secant of $A$ contains exactly one black point of $\mathcal B\cap\Sigma$.

Fix $e\in A$. Each $e'\in A\setminus\{e\}$ determines a unique point of $\mathcal B\cap\Sigma$ on the line $ee'$. Conversely, for each $b\in\mathcal B\cap\Sigma$, equation~\eqref{eq:linebalance} shows that $eb$ contains exactly one further point of $A$. These correspondences give a bijection between $A\setminus\{e\}$ and $\mathcal B\cap\Sigma$. Hence $|\mathcal B\cap\Sigma|=|A|-1$.

The set $\mathcal B\cap\Sigma$ is disjoint from $A$ and meets every secant of $A$, so $A$ is generalized hyperfocused. Theorem~\ref{thm:hyperfocused} now gives $|A|\in\{2,4\}$.

Fix a dual line $\ell$ through two points of $\mathcal{W}$. By assumption, $|\mathcal{W}\cap\ell|=2$. Each of the $p+1$ planes through $\ell$ contains at most four points of $\mathcal{W}$, and hence at most two outside $\ell$. Every point of $\mathcal{W}\setminus\ell$ lies in exactly one of these planes. Therefore $|\mathcal{W}|\le2+2(p+1)=2p+4<1+p(p-1)/2$ for $p\ge11$, contrary to $|\mathcal{W}|=1+p(p-1)/2$.
\end{proof}

\Needspace{6\baselineskip}
\begin{proposition}\label{prop:noquadric}
The set $\mathcal{W}$ is not contained in any quadric defined over $K$.
\end{proposition}

\begin{proof}
Suppose that $\mathcal{W}\subseteq V(Q)$ for some nonzero homogeneous quadratic form $Q$ over $K$. We first show that $V(Q)$ contains no black point, that is, $\mathcal B\cap V(Q)=\varnothing$.

Suppose that $b\in\mathcal B\cap V(Q)$, and choose rational coordinates with $b=(0:0:0:1)$. Since $Q(b)=0$, we can write
\[
Q(X,Y,Z,T)=TL_1(X,Y,Z)+L_2(X,Y,Z),
\]
where $L_1\in\mathscr S_1$ and $L_2\in\mathscr S_2$ are not both zero.

For each $P_i\in D$, equation~\eqref{eq:linebalance} shows that the line $\langle b,P_i\rangle$ contains at least two points of $\mathcal{W}$. Choose two distinct points from this intersection. Using the representative $a_i$, write them as $[a_i:t_{i1}]$ and $[a_i:t_{i2}]$, where $t_{i1}\ne t_{i2}$. Since $Q(a_i,t_{i1})=Q(a_i,t_{i2})=0$, we obtain
\begin{equation}\label{eq:quadric-line-values}
\begin{aligned}
t_{i1}L_1(a_i)+L_2(a_i)=0, \quad
t_{i2}L_1(a_i)+L_2(a_i)=0.
\end{aligned}
\end{equation}
Since $t_{i1}\ne t_{i2}$, a simple calculation using \eqref{eq:quadric-line-values} gives $L_1(a_i)=L_2(a_i)=0$.

Thus $D\subseteq V(L_1)\cap V(L_2)$. Since $(L_1,L_2)\ne(0,0)$, $D$ lies on a curve of degree at most two. This contradicts Proposition~\ref{prop:window}, which states that every curve containing $D$ has degree at least $\lceil(p+1)/3\rceil\ge4$. Therefore $\mathcal B\cap V(Q)=\varnothing$.

By Lemma~\ref{lem:nocap}, there exists a dual line $\ell$ with $|\mathcal{W}\cap\ell|\ge3$. The restriction $Q|_\ell$ is a homogeneous quadratic in two variables that vanishes at three distinct points. Hence $Q|_\ell=0$, so $\ell\subseteq V(Q)$. On the other hand, \eqref{eq:linebalance} gives $\mathcal B\cap\ell\ne\varnothing$. Since $\ell\subseteq V(Q)$, this contradicts $\mathcal B\cap V(Q)=\varnothing$.
\end{proof}

\Needspace{9\baselineskip}
\subsection{Construction of a quadric from the moment identities}

We first prove an independence statement for the projected points whose lines contain at least three points of $\mathcal{W}$.

\begin{lemma}\label{lem:exceptionindependence}
For every $b\in\mathcal B$, the set \(D_{\ge3}=\{P_i\in D:r_i\ge3\}\) is independent in degree $(p-5)/2$; equivalently, $\mathcal R_{(p-5)/2}(D_{\ge3})=\{0\}$.

\end{lemma}

\begin{proof}
Let $s$ be the least degree of a curve containing $D$, and choose $F$ as in Lemma~\ref{lem:window-components}, with $\deg F=s$ and $D\subseteq V(F)$. Put $t_0=(p-2s+1)/2$. By \eqref{eq:window}, we have $t_0\ge2$, and every irreducible factor of $F$ has degree at least $p-2s+1=2t_0$, so $s\ge2t_0$. Since $p=2s+2t_0-1$, equations \eqref{eq:excess} and \eqref{eq:window} give
\begin{equation}\label{eq:exceptionbound}
\begin{aligned}
|D_{\ge3}|\le 1+\frac{p(p-1)}2-2k\le 2t_0^2-3t_0-s+2<st_0,
\end{aligned}
\end{equation}
where the last inequality follows from $s\ge2t_0$ and $t_0\ge2$.

Suppose, for a contradiction, that $0\ne c\in\mathcal R_{(p-5)/2}(D_{\ge3})$, and write $Y=\{P\in D_{\ge3}:c_P\ne0\}$. Since the coefficients outside $Y$ are zero, $c|_Y$ is a nonzero relation in $\mathcal R_{(p-5)/2}(Y)$ whose coordinates are all nonzero.

We first show that no curve of degree less than $t_0$ contains $Y$. Suppose that a nonzero form $G$ with $\deg G<t_0$ vanishes on $Y$. Since every irreducible factor of $F$ has degree at least $2t_0>\deg G$, we have $\gcd(F,G)=1$. Moreover,
\[
\frac{p-5}{2}=s+t_0-3\ge s+\deg G-2.
\]
Corollary~\ref{cor:ciindependence} therefore applies to $Y\subseteq V(F,G)$ and gives $\mathcal R_{(p-5)/2}(Y)=\{0\}$, contrary to $c|_Y\ne0$.

The least degree of a curve containing $Y$ is therefore at least $t_0$. Applying Proposition~\ref{prop:relationbound} to $c|_Y$ with $d=(p-5)/2$ and $t=t_0$, we obtain
\[
|Y|\ge t_0\left(\frac{p+1}{2}-t_0\right)=st_0.
\]
This contradicts $|Y|\le|D_{\ge3}|<st_0$ from \eqref{eq:exceptionbound}, and proves the assertion.
\end{proof}

\begin{proposition}\label{prop:liftquadric}
Let $b\in\mathcal B$, with the coordinates and projection notation as above. Suppose that there exist forms $L\in\mathscr S_1$ and $B\in\mathscr S_2$ such that the moment vectors satisfy
\[
m_2=Lm_1,\qquad m_3=Bm_1
\quad\text{on }D.
\]
Then there is a quadric over $K$ containing all the points of $\mathcal{W}$.
\end{proposition}

\begin{proof}
Let
\begin{equation}\label{eq:liftquadric}
Q(X,Y,Z,T)=T^2-L(X,Y,Z)T+L(X,Y,Z)^2-B(X,Y,Z).
\end{equation}
For each $i$ with $r_i=2$, equation~\eqref{eq:twopoint} and the assumed moment identities give
\[
L(a_i)=t_{i1}+t_{i2},\qquad
B(a_i)=t_{i1}^2+t_{i1}t_{i2}+t_{i2}^2.
\]
Consequently, $Q(a_i,T)=(T-t_{i1})(T-t_{i2})$, and hence $\mathcal{W}\cap\langle b,P_i\rangle\subseteq V(Q)$.

Suppose that $\mathcal{W}'=\mathcal{W}\setminus V(Q)\ne\varnothing$. Since $(c_{i\nu})\in\mathcal R_{p-2}(\mathcal{W})$ by \eqref{eq:fullrelation}, the multiplication rule \eqref{eq:relationmultiply} gives
\[
\bigl(c_{i\nu}Q(a_i,t_{i\nu})\bigr)
\in\mathcal R_{p-4}(\mathcal{W}).
\]
These coefficients vanish on $\mathcal{W}\setminus \mathcal{W}'$ and are nonzero on $\mathcal{W}'$, since every $c_{i\nu}\ne0$. Thus, for every $A\in K[X,Y,Z,T]_{p-4}$, we have
\begin{equation}\label{eq:quadric-relation}
\sum_{\substack{i,\nu\\ [a_i:t_{i\nu}]\in \mathcal{W}'}}
c_{i\nu}Q(a_i,t_{i\nu})A(a_i,t_{i\nu})=0.
\end{equation}

Since $\mathcal{W}\cap\langle b,P_i\rangle\subseteq V(Q)$ whenever $r_i=2$, no point of $\mathcal{W}'$ lies on such a line. As $r_i\ge2$ for every $i$, each point of $\mathcal{W}'$ projects to a point $P_i$ with $r_i\ge3$. Hence $\pi_{b,H}(\mathcal{W}')\subseteq D_{\ge3}$. By Lemma~\ref{lem:exceptionindependence}, $D_{\ge3}$ is independent in degree $(p-5)/2$, and hence in degree $v$ for every integer $v\ge (p-5)/2$.

Fix $P_i\in\pi_{b,H}(\mathcal{W}')$, and write
\[
\mathcal{W}'\cap\langle b,P_i\rangle
=\{[a_i:t_\nu]:1\le\nu\le r\},
\]
where $t_1,\ldots,t_r$ are pairwise distinct. Denote the coefficients in \eqref{eq:quadric-relation} corresponding to these points by $d_1,\ldots,d_r$; each is nonzero. By Lemma~\ref{lem:linecap}, we have $r\le(p-1)/2$, and hence
\[
p-4-j\ge p-3-\frac{p-1}{2}=\frac{p-5}{2}
\qquad(0\le j<r).
\]
For each such $j$, the independence of $D_{\ge3}$ in degree $p-4-j$ implies that there exists $G_j\in\mathscr S_{p-4-j}$ satisfying
\begin{equation}\label{eq:interpolation-values}
G_j(a_i)=1,\qquad
G_j(a_\ell)=0
\quad\text{for all }P_\ell\in D_{\ge3}\setminus\{P_i\}.
\end{equation}

Since $\deg(T^jG_j)=p-4$, substituting $A=T^jG_j$ in \eqref{eq:quadric-relation} and using \eqref{eq:interpolation-values} gives
\[
\sum_{\nu=1}^r d_\nu t_\nu^{\,j}=0
\qquad(0\le j<r).
\]
The coefficient matrix is a Vandermonde matrix with determinant $\prod_{1\le\mu<\nu\le r}(t_\nu-t_\mu)\ne0$, since the $t_\nu$ are pairwise distinct. Thus $d_1=\cdots=d_r=0$, a contradiction. Hence $\mathcal{W}'=\varnothing$ and $\mathcal{W}\subseteq V(Q)$.

\end{proof}

\Needspace{12\baselineskip}

We prove Theorem~\ref{thm:main} in its equivalent projective form, Theorem~\ref{thm:main}$^\ast$.

\begin{proof}[Proof of Theorem~\ref{thm:main}]
Every set satisfying \eqref{eq:hyp} is a cylinder for $p=2,3$ by \cite{BDMS}, and for $p=5,7$ by \cite[Propositions~41 and~52]{KM}, respectively. It remains to consider $p\ge11$. Suppose, for a contradiction, that $S$ satisfies \eqref{eq:hyp} but is not a cylinder. Lemma~\ref{lem:incidence} gives $\mathcal B\ne\varnothing$, and Proposition~\ref{prop:simple} gives $|\mathcal{W}|=1+p(p-1)/2$ and $w_b=1$ for every $b\in\mathcal B$.

Choose $b\in\mathcal B$, with the projection notation as above. By Theorem~\ref{thm:globalmoments}, there exist forms $L\in\mathscr S_1$ and $B\in\mathscr S_2$ such that
\[
m_2=Lm_1,\qquad m_3=Bm_1
\quad\text{on }D.
\]
Proposition~\ref{prop:liftquadric} therefore shows that $\mathcal{W}$ is contained in a quadric over $K$, contradicting Proposition~\ref{prop:noquadric}. Hence $S$ is a cylinder.
\end{proof}



\section{Proof of R\'edei's factorization conjecture}\label{sec:redei}

For a nontrivial normalized factorization $A\oplus B=\Fp^3$, we may assume that $|A|=p$ and $|B|=p^2$. If $A$ spans $\Fp^3$, we prove that $B$ is $p$-divisible; Theorem~\ref{thm:main} then implies that $B$ is a cylinder. The remaining argument applies R\'edei's planar factorization theorem in a two-dimensional quotient of $\Fp^3$ to prove Theorem~\ref{thm:redei}.

For column vectors $u=(u_1,\ldots,u_m)^{\mathsf T}$ and $v=(v_1,\ldots,v_m)^{\mathsf T}$ over $\Fp$, write
\[
u\cdot v=\sum_{i=1}^m u_iv_i.
\]
The same notation applies when the entries of either vector are polynomials over $\Fp$.

\subsection{Factorizations and power sums of linear forms}

The Fourier characterization of tilings and the equidistribution criterion in \cite[Lemma~1.2 and Proposition~2.2]{HIPRV} give the following lemma; see also \cite[Corollary~5.1.6]{MattheusThesis}.

\begin{lemma}\label{lem:factorfourier}
Suppose that $A\oplus B=\Fp^m$, where $m\ge1$. For every $0\ne u\in\Fp^m$, there exists $C\in\{A,B\}$ such that
\[
|\{c\in C:u\cdot c=t\}|=\frac{|C|}{p}
\qquad(t\in\Fp).
\]
\end{lemma}

For $C\subseteq\Fp^3$ and $X=(X_1,X_2,X_3)^{\mathsf T}$, define
\[
M_{C,0}=\sum_{c\in C}1\in\Fp,
\qquad
M_{C,j}(X)=\sum_{c\in C}(c\cdot X)^j
\in\Fp[X_1,X_2,X_3]\quad(j\ge1).
\]

\begin{lemma}\label{lem:firstmoment}
Let $A\subseteq\Fp^3$ satisfy $|A|=p$ and $0\in A$. If $A$ spans $\Fp^3$, then $M_{A,r}\ne0$ for some $1\le r\le p-2$.
\end{lemma}

\begin{proof}
Since $0\in A$ and $\langle A\rangle=\Fp^3$, we have $|A\setminus\{0\}|\ge3$ and hence $p\ge5$. Choose $a\in A\setminus\{0\}$ and $v\in\Fp^3\setminus\{0\}$ with $v\cdot a=0$. The spanning assumption gives $A\nsubseteq\{x\in\Fp^3:v\cdot x=0\}$, whereas $v\cdot0=v\cdot a=0$. Therefore
\begin{equation}\label{eq:firstmoment-values}
2\le|\{v\cdot x:x\in A\}|\le p-1.
\end{equation}
Suppose that
\[
\sum_{x\in A}(v\cdot x)^j=0
\qquad\left(1\le j\le\frac{p-1}{2}\right).
\]
A consequence of R\'edei's theorem \cite[Satz~1]{RedeiGaps1947} (see \cite[Section~4.2]{LS} for the power-sum formulation) is that
\[
|\{v\cdot x:x\in A\}|\in\{1,p\},
\]
contrary to \eqref{eq:firstmoment-values}. Consequently,
\[
M_{A,r}(v)=\sum_{x\in A}(v\cdot x)^r\ne0
\qquad\text{for some }1\le r\le\frac{p-1}{2}\le p-2,
\]
which implies $M_{A,r}\ne0$.
\end{proof}

\subsection{Divisibility of the larger factor and proof of R\'edei's conjecture}

\begin{lemma}\label{lem:momentvector}
Let $A\subseteq\Fp^3$ satisfy $0\in A$ and $\langle A\rangle=\Fp^3$. Suppose that $1\le r<p$ and
\[
M_{A,r}\ne0,\qquad M_{A,j}=0\quad(0\le j<r).
\]
Let $X=(X_1,X_2,X_3)^{\mathsf T}$ be a vector of variables and $U=(U_1,U_2,U_3)^{\mathsf T}$, where each $U_i\in\Fp[X_1,X_2,X_3]$ is homogeneous of degree $e$, with $0\le e\le r$. If \((a-b)\cdot X\mid(a-b)\cdot U \) for all distinct \(a,b\in A,\) then $U=qX$ for some $q\in\Fp[X_1,X_2,X_3]$.
\end{lemma}

\begin{proof}
For $a\in A$, put $\ell_a=a\cdot X$. If $a\ne b$, then $\ell_a-\ell_b=(a-b)\cdot X\ne0$. Thus the polynomials $\ell_a$, $a\in A$, are pairwise distinct, and the Lagrange interpolation formula defines
\begin{equation}\label{eq:factorlagrange}
P(T)=\sum_{a\in A}(a\cdot U)
\prod_{b\in A\setminus\{a\}}
\frac{T-\ell_b}{\ell_a-\ell_b}
\in\Fp(X_1,X_2,X_3)[T].
\end{equation}
It follows directly that
\begin{equation}\label{eq:factorinterpolation}
\deg_TP<|A|,\qquad P(\ell_a)=a\cdot U\quad(a\in A).
\end{equation}
The difference of two polynomials satisfying \eqref{eq:factorinterpolation} has degree less than $|A|$ and has the $|A|$ distinct roots $\ell_a$, $a\in A$; hence $P$ is unique.

We prove that every coefficient of $P$ is a polynomial. Write an arbitrary coefficient as $f/g$, where $f,g\in\Fp[X_1,X_2,X_3]$, $g\ne0$, and $\gcd(f,g)=1$. By \eqref{eq:factorlagrange}, every irreducible factor of $g$ is of the form $\lambda(\ell_a-\ell_b)$, where $a,b\in A$, $a\ne b$, and $\lambda\in\Fp^\times$, because $\ell_a-\ell_b$ is linear. Thus, to prove $g\in\Fp^\times$, it suffices to show that \(\ell_a-\ell_b\nmid g\) for all distinct \(a,b\in A.\)

Fix distinct $a,b\in A$ and put $v=a-b$. By hypothesis, there exists $q_v\in\Fp[X_1,X_2,X_3]$ such that
\begin{equation}\label{eq:factorlinearquotient}
v\cdot U=q_v(v\cdot X).
\end{equation}

The distinct sets $A\cap(c+\Fp v)$, $c\in A$, partition $A$. We may therefore write $A=\{a_1,\ldots,a_{|A|}\}$, with the indices chosen so that
\begin{equation}\label{eq:factorenumeration}
a_j-a_i\in\Fp v
\quad\Longrightarrow\quad
a_k-a_i\in\Fp v
\qquad(1\le i\le k\le j\le |A|).
\end{equation}

Define $\Delta_{i,j}\in\Fp(X_1,X_2,X_3)$ by
\begin{equation}\label{eq:factordifference}
\begin{aligned}
\Delta_{i,i}&=a_i\cdot U
&& (1\le i\le |A|),\\
\Delta_{i,j}&=\frac{\Delta_{i+1,j}-\Delta_{i,j-1}}
{\ell_{a_j}-\ell_{a_i}}
&& (1\le i<j\le |A|).
\end{aligned}
\end{equation}
The polynomial $P$ in \eqref{eq:factorinterpolation} can also be written as
\begin{equation}\label{eq:factornewton}
P(T)=\sum_{j=1}^{|A|}\Delta_{1,j}
\prod_{i=1}^{j-1}(T-\ell_{a_i}),
\end{equation}
by the Newton interpolation formula \cite[ Proposition~32]{deBoor}. Since $\ell_{a_i}-\ell_{a_j}\ne0$ for $i\ne j$, all divisions in \eqref{eq:factordifference} are valid in $\Fp(X_1,X_2,X_3)$, and the same interpolation formula applies over this field.

Suppose that $a_j-a_i\in\Fp v$ with $i<j$. By \eqref{eq:factorenumeration}, write $a_k=a_i+t_kv$ with $t_k\in\Fp$ for $i\le k\le j$. Equation~\eqref{eq:factorlinearquotient} then gives
\begin{equation}\label{eq:factorlinevalues}
\ell_{a_k}=\ell_{a_i}+t_k(v\cdot X),\qquad
a_k\cdot U=a_i\cdot U+t_kq_v(v\cdot X)
\quad(i\le k\le j).
\end{equation}
Since $t_i,\ldots,t_j$ are pairwise distinct, \eqref{eq:factorlinevalues} and \eqref{eq:factordifference} imply
\begin{equation}\label{eq:factorlinedifferences}
\Delta_{i,j}=
\begin{cases}
q_v,&j=i+1,\\
0,&j\ge i+2,
\end{cases}
\qquad\text{if }a_j-a_i\in\Fp v.
\end{equation}
Indeed, for $i\le k<j$,
\[
\Delta_{k,k+1}
=\frac{(t_{k+1}-t_k)q_v(v\cdot X)}
{(t_{k+1}-t_k)(v\cdot X)}
=q_v.
\]
Substitution in \eqref{eq:factordifference} gives $\Delta_{k,k+2}=0$, and induction in the same formula gives $\Delta_{i,j}=0$ for $j\ge i+2$.

We prove by induction on $j-i$ that there exist $f_{i,j},g_{i,j}\in\Fp[X_1,X_2,X_3]$ such that
\[
\Delta_{i,j}=\frac{f_{i,j}}{g_{i,j}},\qquad
g_{i,j}\ne0,\qquad
v\cdot X\nmid g_{i,j}.
\]
For $i=j$, take $f_{i,i}=a_i\cdot U$ and $g_{i,i}=1$, as permitted by \eqref{eq:factordifference}. If $i<j$ and $a_j-a_i\in\Fp v$, then \eqref{eq:factorlinedifferences} gives $\Delta_{i,j}\in\{q_v,0\}$, so we may take $f_{i,j}=\Delta_{i,j}$ and $g_{i,j}=1$.

Now suppose that $i<j$ and $a_j-a_i\notin\Fp v$. The linear forms $v\cdot X$ and $\ell_{a_j}-\ell_{a_i}=(a_j-a_i)\cdot X$ are not proportional, so
\[
v\cdot X\nmid\ell_{a_j}-\ell_{a_i}.
\]
By the induction hypothesis,
\[
g_{i+1,j}g_{i,j-1}\ne0,\qquad
v\cdot X\nmid g_{i+1,j}g_{i,j-1}.
\]
Put
\[
\begin{aligned}
f_{i,j}
&=f_{i+1,j}g_{i,j-1}-f_{i,j-1}g_{i+1,j},\\
g_{i,j}
&=g_{i+1,j}g_{i,j-1}(\ell_{a_j}-\ell_{a_i}).
\end{aligned}
\]
Equation~\eqref{eq:factordifference} gives $\Delta_{i,j}=f_{i,j}/g_{i,j}$. Since $v\cdot X$ is irreducible, the preceding divisibility conditions imply
\[
g_{i,j}\ne0,\qquad v\cdot X\nmid g_{i,j},
\]
which completes the induction.

Since the products in \eqref{eq:factornewton} have polynomial coefficients, we may write
\[
\frac{f}{g}=\frac{\widetilde f}{\widetilde g},
\qquad
\widetilde g=\prod_{j=1}^{|A|}g_{1,j},
\]
with $\widetilde f\in\Fp[X_1,X_2,X_3]$. The induction above gives $\widetilde g\ne0$ and $v\cdot X\nmid\widetilde g$. Since $f\widetilde g=\widetilde f g$ and $\gcd(f,g)=1$, we have $g\mid\widetilde g$ and hence $v\cdot X\nmid g$. Since $v\cdot X=\ell_a-\ell_b$ and the distinct points $a,b\in A$ were arbitrary, $g$ has no feasible irreducible factrors. Therefore $g\in\Fp^\times$ and $P\in\Fp[X_1,X_2,X_3][T]$.

Write $P(T)=\sum_{j=0}^{|A|-1}P_jT^j$. By \eqref{eq:factorlagrange}, each nonzero $P_j$ is homogeneous of degree \(e+(|A|-1-j)-(|A|-1)=e-j.\) Since $P_j\in\Fp[X_1,X_2,X_3]$, we have $P_j=0$ for $j>e$. Moreover, $0\in A$ and \eqref{eq:factorinterpolation} give $P_0=P(0)=0$.

Suppose that $\deg_TP=j\ge2$. Then $2\le j\le e\le r$, so $1\le r-j+1<r$ and $M_{A,r-j+1}=0$. Differentiating $M_{A,r-j+1}=0$ with respect to $X_1,X_2,X_3$ and using \eqref{eq:factorinterpolation}, we obtain
\[
\begin{aligned}
0
=\sum_{\ell=1}^{3}U_\ell
  \frac{\partial M_{A,r-j+1}}{\partial X_\ell}
=(r-j+1)\sum_{a\in A}(a\cdot U)\ell_a^{\,r-j}&=(r-j+1)\sum_{i=0}^{j}P_iM_{A,r-j+i}.
\end{aligned}
\]
For $0\le i<j$, we have $0\le r-j+i<r$ and therefore $M_{A,r-j+i}=0$. Hence \((r-j+1)P_jM_{A,r}=0.\) Since $1\le r-j+1<p$, we conclude that $P_jM_{A,r}=0$, contrary to $P_j\ne0$ and $M_{A,r}\ne0$. We have therefore shown that $\deg_TP\le1$. Since $P_0=0$, we have $P(T)=qT$ with $q=P_1$. Substituting into \eqref{eq:factorinterpolation} gives
\[
a\cdot(U-qX)=0\qquad(a\in A).
\]
As $\langle A\rangle=\Fp^3$, we conclude that $U=qX$.
\end{proof}

\begin{proposition}\label{prop:factorpdivisible}
Let $A\oplus B=\Fp^3$ with $|A|=p$ and $0\in A$. If $A$ spans $\Fp^3$, then $B$ is $p$-divisible.
\end{proposition}

\begin{proof}
The factorization gives $|A||B|=p^3$, so
\begin{equation}\label{eq:factor-zero-moments}
|A|=p,\qquad |B|=p^2,\qquad M_{A,0}=M_{B,0}=0.
\end{equation}
Let $r$ be the least index for which $M_{A,r}\ne0$. By Lemma~\ref{lem:firstmoment} and \eqref{eq:factor-zero-moments},
\begin{equation}\label{eq:factor-first-moment}
1\le r\le p-2,\qquad M_{A,j}=0\quad(0\le j<r).
\end{equation}
Suppose that $B$ is not $p$-divisible. Then there exist $0\ne u_0\in\Fp^3$ and $t_*\in\Fp$ such that
\[
|\{b\in B:u_0\cdot b=t_*\}|\not\equiv0\pmod p.
\]
The equality $z^{p-1}=1$ for $z\in\Fp^\times$ gives
\begin{equation}\label{eq:factor-plane-function}
1-(u_0\cdot x-t_*)^{p-1}
=\begin{cases}
1,&u_0\cdot x=t_*,\\
0,&u_0\cdot x\ne t_*,
\end{cases}
\qquad(x\in\Fp^3).
\end{equation}
Summing \eqref{eq:factor-plane-function} over $B$ in $\Fp$ and using \eqref{eq:factor-zero-moments}, we obtain
\begin{equation}\label{eq:factor-nonzero-moment}
\begin{aligned}
0&\ne\sum_{b\in B}\bigl(1-(u_0\cdot b-t_*)^{p-1}\bigr)\\
&=-\sum_{j=1}^{p-1}\binom{p-1}{j}(-t_*)^{p-1-j}M_{B,j}(u_0).
\end{aligned}
\end{equation}
Equation~\eqref{eq:factor-nonzero-moment} gives $M_{B,d}(u_0)\ne0$ for some $1\le d\le p-1$. Put $F=M_{A,r}$ and $G=M_{B,d}$, so that
\[
\deg F=r,\qquad\deg G=d,\qquad FG\ne0.
\]

Fix distinct $a,a'\in A$. For every $0\ne u\in\Fp^3$ with $u\cdot(a-a')=0$, we have
\[
u\cdot a=u\cdot a',\qquad
|\{x\in A:u\cdot x=u\cdot a\}|
\ge2>\frac{|A|}{p}=1.
\]
Lemma~\ref{lem:factorfourier} therefore gives
\begin{equation}\label{eq:factor-uniform-B}
|\{b\in B:u\cdot b=t\}|=p\quad(t\in\Fp),
\qquad
G(u)=p\sum_{t\in\Fp}t^d=0.
\end{equation}

By \eqref{eq:factor-uniform-B}, $G$ vanishes at all $p+1$ points of the line $(a-a')\cdot X=0$ in $\PG(2,p)$. Since $\deg G=d<p+1$, its restriction to this line is identically zero. Therefore
\begin{equation}\label{eq:factorline-divisibility}
(a-a')\cdot X\mid G
\qquad\text{for all distinct }a,a'\in A.
\end{equation}

For any $0\ne u\in\Fp^3$, Lemma~\ref{lem:factorfourier} gives
\[
F(u)=\sum_{t\in\Fp}t^r=0
\quad\text{or}\quad
G(u)=p\sum_{t\in\Fp}t^d=0,
\]
where $\sum_{t\in\Fp}t^r=0$ since $1\le r\le p-2$. Together with $F(0)=G(0)=0$, this gives
\begin{equation}\label{eq:factorpointwise}
F(u)G(u)=0\qquad(u\in\Fp^3).
\end{equation}

By \cite[Lemma~7.2.4]{BFS} (see \cite[Section~2]{MR} for a proof), the homogeneous vanishing ideal of $\PG(2,p)$ over $\Fp$ is
\begin{equation}\label{eq:factor-projective-ideal}
\bigl(X_i^pX_j-X_iX_j^p:1\le i<j\le3\bigr)
\subseteq\Fp[X_1,X_2,X_3].
\end{equation}
By \eqref{eq:factorpointwise}, $FG$ belongs to the ideal in \eqref{eq:factor-projective-ideal}. Each of these forms has degree $p+1$, so $FG\ne0$ implies $r+d=\deg(FG)\ge p+1$. Since $FG$ is homogeneous of degree $r+d$, there exist homogeneous polynomials $U_1,U_2,U_3$ of degree $r+d-p-1$ such that
\begin{equation}\label{eq:factordeterminant}
\begin{aligned}
FG
&=(X_2^pX_3-X_2X_3^p)U_1
 +(X_3^pX_1-X_3X_1^p)U_2\\
&\quad +(X_1^pX_2-X_1X_2^p)U_3\\
&=\det(X^{[p]},X,U),
\end{aligned}
\end{equation}
where $X^{[p]}=(X_1^p,X_2^p,X_3^p)^{\mathsf T}$ and $U=(U_1,U_2,U_3)^{\mathsf T}$. Since $d=\deg G\le p-1$,
\begin{equation}\label{eq:factor-coefficient-degree}
0\le r+d-p-1\le r-2.
\end{equation}
In particular, $r\ge2$.

Fix distinct $a,a'\in A$, and choose a basis $u,w$ of
\[
\{x\in\Fp^3:(a-a')\cdot x=0\}.
\]
Substitute $X=Yu+Zw$, where $Y,Z$ are variables. Since $u,w\in\Fp^3$, we have $X^{[p]}=Y^pu+Z^pw$. Equations \eqref{eq:factorline-divisibility} and \eqref{eq:factordeterminant} give
\[
0=(FG)(Yu+Zw)
 =(Y^pZ-Z^pY)\det\bigl(u,w,U(Yu+Zw)\bigr).
\]
Since $Y^pZ-Z^pY\ne0$ in $\Fp[Y,Z]$, we obtain $\det(u,w,U(Yu+Zw))=0$. The nonzero linear forms $\det(u,w,X)$ and $(a-a')\cdot X$ have the same kernel $\langle u,w\rangle$ and are therefore proportional. It follows that
\[
(a-a')\cdot U(Yu+Zw)=0
\qquad\text{in }\Fp[Y,Z].
\]
Hence
\begin{equation}\label{eq:factor-vector-divisibility}
(a-a')\cdot X\mid(a-a')\cdot U
\qquad\text{for all distinct }a,a'\in A.
\end{equation}

By \eqref{eq:factor-first-moment}, \eqref{eq:factor-coefficient-degree}, and \eqref{eq:factor-vector-divisibility}, Lemma~\ref{lem:momentvector} gives $U=qX$ for some $q\in\Fp[X_1,X_2,X_3]$. Substituting into \eqref{eq:factordeterminant}, we obtain
\[
FG=\det(X^{[p]},X,qX)=0,
\]
contrary to $FG\ne0$. Therefore $B$ is $p$-divisible.
\end{proof}

We now prove Theorem~\ref{thm:redei}.

\begin{proof}[Proof of Theorem~\ref{thm:redei}]
The normalization gives $0\in A\cap B$. If $|A|=1$ or $|B|=1$, then $A=\{0\}$ or $B=\{0\}$, respectively, and the assertion follows. Assume that $|A|,|B|>1$. Since $|A||B|=p^3$, interchanging $A$ and $B$ if necessary gives
\[
|A|=p,\qquad |B|=p^2.
\]
If $\dim\langle A\rangle\le2$, then $A$ is contained in a plane through zero. We may therefore assume that $\langle A\rangle=\Fp^3$. Proposition~\ref{prop:factorpdivisible} and Theorem~\ref{thm:main} then imply that $B$ is a cylinder. Thus there exists a subspace $L\le\Fp^3$, with $\dim L=1$, such that $B$ is the union of $p$ distinct cosets of $L$.

Following the quotient argument of Szab\'o and Ward \cite{SzaboWard}, consider the linear map $\pi_L:\Fp^3\to\Fp^3/L$. Since $B$ is a union of cosets of $L$,
\[
B+L=B,\qquad B=\pi_L^{-1}(\pi_L(B)),\qquad |\pi_L(B)|=p.
\]
If $\pi_L(a)=\pi_L(a')$ for $a,a'\in A$, then $a-a'\in L$ and, for every $b\in B$,
\[
a+b=a'+\bigl(b+(a-a')\bigr),\qquad b+(a-a')\in B.
\]
Uniqueness in $A\oplus B=\Fp^3$ implies $a=a'$; hence $|\pi_L(A)|=|A|=p$. Since $A+B=\Fp^3$, the addition map $\pi_L(A)\times\pi_L(B)\to\Fp^3/L$ is surjective. The equality
\[
|\pi_L(A)\times\pi_L(B)|=p^2=|\Fp^3/L|
\]
therefore implies that this map is bijective, so
\[
\pi_L(A)\oplus\pi_L(B)=\Fp^3/L.
\]

We also have $0=\pi_L(0)\in\pi_L(A)\cap\pi_L(B)$. R\'edei's planar factorization theorem \cite[Satz~4]{RedeiGaps1947} (see also \cite[Theorem~23]{RedeiBookEnglish}) therefore implies that at least one of $\pi_L(A)$ and $\pi_L(B)$ is a one-dimensional subspace of $\Fp^3/L$. Since $\langle A\rangle=\Fp^3$ and $\pi_L$ is linear and surjective, we have
\[
\langle\pi_L(A)\rangle=\Fp^3/L,\qquad \dim(\Fp^3/L)=2.
\]
Thus, $\pi_L(A)$ is not a one-dimensional subspace and $\pi_L(B)$ is a one-dimensional subspace. As $\ker\pi_L=L$, its inverse image $B=\pi_L^{-1}(\pi_L(B))$ is a subspace with
\[
\dim B=\dim L+\dim\pi_L(B)=2.
\]
Hence, $B$ is a plane through zero.
\end{proof}

\begin{remark}
A normalized factorization is called a \emph{full rank tiling} if both factors span the whole space. For every prime $p\ge5$, Dinitz \cite[Proposition~12]{Dinitz} constructed a full rank tiling of $\Fp^4$. In dimension one, $|A||B|=p$ implies $A=\{0\}$ or $B=\{0\}$, since $0\in A\cap B$. In dimension two, either one factor is $\{0\}$ or $|A|=|B|=p$; in the latter case, the planar result recalled in Section~\ref{sec:intro} implies that one factor is a line. Theorem~\ref{thm:redei} excludes full rank tilings in dimension three. Hence, for every prime $p\ge5$, four is the least dimension in which a full rank tiling exists.
\end{remark}

\section*{Acknowledgments}

We are grateful to Yue Zhou for pointing out the connection between the strong cylinder conjecture and R\'edei's factorization conjecture, an observation that led us to a proof of the latter.

This research was supported by the National Key Research and Development Program of China (Grant No.~2025YFA1017700). The National Natural Science Foundation of China supported Haihua Deng (Grant No.~123B2011) and Weicong Li (Grant No.~12301422). Haihua Deng also received support from the National Postdoctoral Program for Innovative Talents (Grant No.~BX20250059).

\section*{AI Declaration}

During the preparation of this work, the authors used OpenAI GPT-6 Astra to help develop the proofs, particularly for the interpolation arguments in Section \ref{sec:ci}, write the Lean 4 formalizations \footnote{https://github.com/WeicongLi13/weicongli.github.com/tree/main/cylinder-lean}, and prepare the manuscript. The authors take full responsibility for the contents of the paper.


\begingroup
\hypersetup{citecolor=black,urlcolor=black,linkcolor=black}
\color{black}
\renewcommand*{\backref}[1]{}

\endgroup

\end{document}